\documentclass[11pt]{article}
\usepackage{latexsym}
\usepackage{amsfonts}
\usepackage{amsmath}
\usepackage{epsfig}

\newtheorem{thm}{Theorem} \newtheorem{lemma}[thm]{Lemma}
\newtheorem{prop}[thm]{Proposition} \newtheorem{cor}[thm]{Corollary}

\newcommand{\prob}{\mbox{\bf P}}
\newcommand{\E}{{\bf E\,}}

\newcommand{\diam}{\mbox{\rm diam}}

\newcommand{\eps}{\varepsilon}
\newcommand{\be}{\begin{equation}}
\newcommand{\ee}{\end{equation}}
\newcommand{\C}{{\cal{C}}}
\newcommand{\G}{{\cal{G}}}

\newcommand{\F}{{\cal{F}}}
\newcommand{\A}{{\cal{A}}}
\renewcommand{\L}{{\cal{L}}}
\renewcommand{\S}{{\cal{S}}}
\newcommand{\B}{{\cal{B}}}

\newcommand{\Ne}{ {\bf{N}} }
\newcommand{\Net}{ {\bf{ \widetilde{N}}} }
\newcommand{\NN} { { \widetilde{N}} }

\newcommand{\Ac}{ {\bf{A}} }
\newcommand{\Act} { {\bf { \widetilde{A}}} }
\renewcommand{\AA} { {\widetilde{A}} }

\newcommand{\Simp}{{\cal S}imple}
\newcommand{\Tm}{T_{{\rm mix}}}

\def\qed{\relax\ifmmode\hskip2em \Box\else\unskip\nobreak\hfill
$\Box$\fi}

\newcounter{mycount}

\title{Critical percolation on random regular graphs}
\author{{\sc Asaf Nachmias and Yuval Peres\thanks{U.C. Berkeley. Research of both authors
supported in part by NSF grants \#DMS-0244479 and
\#DMS-0104073.}}}
\begin{document}

\maketitle

\begin{abstract}
We describe the component sizes in critical independent $p$-bond
percolation on a random $d$-regular graph on $n$ vertices, where
$d\geq 3$ is fixed and $n$ grows. We prove {\em mean-field}
behavior around the critical probability $p_c={1 \over d-1}$.

In particular, we show that there is a scaling window of width
$n^{-1/3}$ around $p_c$ in which the sizes of the largest
components are roughly $n^{2/3}$ and we describe their limiting
joint distribution. We also show that for the subcritical regime,
i.e. $p = (1-\eps(n)) p_c $ where $\eps(n)=o(1)$ but
$\eps(n)n^{1/3} \to \infty$, the sizes of the largest components
are concentrated around an explicit function of $n$ and $\eps(n)$
which is of order $o(n^{2/3})$. In the supercritical regime, i.e.
$p = (1+\eps(n)) p_c$ where $\eps(n)=o(1)$ but $\eps(n)n^{1/3} \to
\infty$, the size of the largest component is concentrated around
the value ${2d \over d-2} \eps(n) n$ and a duality principle
holds: other component sizes are distributed as in the subcritical
regime.

\end{abstract}

\section{Introduction}

Let $d\geq 3$ be a fixed integer, $n>0$ an integer such that $dn$
is even, and $p \in (0,1)$. Let $G(n,d,p)$ be a random graph on
$n$ vertices obtained by drawing uniformly a random $d$-regular
graph on $n$ vertices and then performing independent $p$-bond
percolation on it, i.e., we independently retain each edge with
probability $p$ and delete it with probability $1-p$. Alon,
Benjamini and Stacey proved in \cite{ABS} that the model $G(n,d,{c
\over d-1})$ exhibits a phase transition as $c$ grows: the
cardinality of the largest component $\C_1$ is of order $\log n$
for $c<1$ and of order $n$ for $c>1$.

Recall that similar behavior is exhibited in the random graph
$G(n,p)$, introduced by Erd\H{o}s and R\'enyi~\cite{ER}. They
discovered that as $c$ grows, $G(n,c/n)$ exhibits a {\em double
jump\/}: the cardinality of the largest component $\C_1$ is of
order $\log n$
 for $c<1$,  of order $n^{2/3}$ for $c=1$ and linear in $n$ for
 $c>1$.  In fact, for the critical case  $c=1$ the argument in
\cite{ER} only established the lower bound; the upper bound was
proved much later in \cite{B1}, \cite{L} and \cite{LPW}; see also
\cite{NP} for a simple proof of this upper bound. These works
established the existence of a ``scaling-window'' of width
$n^{-1/3}$ around the point ${1 \over n}$, i.e., for all $p$ of
the form ${1 \over n}(1+O(n^{-1/3}))$ the random variable $|\C_1|
/ n^{2/3}$ converges in distribution to a non-trivial random
variable, and in particular, is not concentrated. Furthermore,
outside of this scaling window, i.e. for $p$ of the form ${1 \over
n} (1 + \eps(n))$ where $\eps(n)=o(1)$ but $\eps(n)n^{1/3} \to
\infty$, the random variable $|\C_1|$ is concentrated around some
known value. This is often called ``mean-field'' behavior around
the critical probability $p_c(n) = {1 \over n}$.

Itai Benjamini (personal communication) asked whether percolation
on a random $d$-regular graph has mean-field behavior. In this
paper we answer his question affirmatively for $d$ fixed and $n$
growing, and give a complete description of the component sizes at
criticality. We establish the existence of a scaling window of
width $n^{-1/3}$ around the critical probability ${1 \over d-1}$
(in which component sizes have a non-trivial limiting
distribution) and show that outside the window the largest
component (and the $\ell$-th largest component as well) is
concentrated. Boris Pittel (personal communication) informed us
that he had obtained similar (but somewhat less
precise) results. \\

Recall (see \cite{B2} and \cite{JLR}) that in the
Erd\H{o}s-R\'enyi random graph $G(n, {1 - \eps(n) \over n})$,
where $\eps(n)>0$ satisfies $\eps(n) \to 0$ and $\eps(n) n^{1/3}
\to \infty$, for any fixed integer $\ell
> 0$ we have
\be \label{gnpsubcrit} { |\C_\ell| \over \psi_n(\eps(n))}
\stackrel{P} \longrightarrow 1 \qquad \hbox{as } n \to \infty \,
,\ee where \be \label{psidef} \psi_n(\eps) = 2 \eps^{-2} \log (
n\eps^{-3} ) \, .\ee

The following proposition provides general upper bounds on the
size of the largest component which are valid for {\em all}
$d$-regular graphs. In particular, part $1$ provides an upper
bound on $|\C_1|$ in the subcritical regime, similar to the one
implied in (\ref{gnpsubcrit}), and part $2$ and $3$ provide upper
bounds for other regimes of $p$.

\begin{prop} \label{dregular} {\bf [General upper bounds]} Let $G$ be a $d$-regular
graph for $d\geq 3$ and denote by $\C_1(G_p)$ the largest
connected component of the random graph obtained by bond
percolation on $G$ with probability $p$. We have
\begin{enumerate}

\item If $p = {1 - \eps(n) \over d-1}$ where $\eps(n)\geq 0$ is a
sequence such that $\eps(n) \to 0$ and $\eps(n) n^{1/3} \to
\infty$, then for any $\eta>0$
$$ \prob \Big (  |\C_1(G_p)| > (1+\eta) {d-2 \over d-1}\psi_n (\eps(n))   \Big ) \to 0 \, ,$$
as $n\to \infty$.

\item If $p \leq {1 \over d-1}$ where $\lambda \leq 0$ then for
any $A>1$
$$ \prob \Big ( |\C_1(G_p)| > An^{2/3} \Big ) \leq {8 \over A^{3/2}} \, .$$

\item There exists a constant $C>0$ such that if $p={1+\eps(n)
\over d-1}$ where $\eps(n) > 0$ then
$$ \E |\C_1(G_p)| \leq C ( n^{2/3} + \eps(n) n ) \, .$$
\end{enumerate}
\end{prop}

For a random regular graph, we can sharpen these upper bounds and
prove corresponding lower bounds. In the following we denote by
$\{ \C _j \}_{j\geq 1}$ the connected components of $G(n,d,p)$
ordered in decreasing size. We emphasize that all the theorems
apply for $d$ fixed and $n$ growing. See Section \ref{secconc} for
further discussion on the case where $d$ grows with $n$.

\begin{thm} \label{criticalupper} {\bf [Critical window bounds]}
Consider $G(n,d,p)$ with $p={1+\lambda n^{-1/3} \over d-1}$ for
some $\lambda \in {\mathbb R}$ where $d\geq 3$ is fixed. Then
there there exist constants $c(\lambda,d)>0$ and
$C(\lambda,d)<\infty$ such that for any $A>0$ and all $n$,

\be \label{part1} \prob ( |\C_1| \geq An^{2/3} ) \leq
{C(\lambda,d) e^{-c(\lambda,d) A^3} \over A} \, .\ee Furthermore,
there exists a constant $D= D(\lambda,d)$ such that for $\delta>0$
small enough and all $n$,

\be \label{part2} \prob \Big ( |\C_1| <\lceil \delta n^{2/3}
\rceil \Big ) \leq D(\lambda,d) \delta^{1/2}  \, .\ee

\end{thm}

The next two theorems describe the largest component behavior
outside of the scaling window. In particular, outside the scaling
window, the largest component is concentrated; however, the
structure of the graphs is quite different depending on whether we
are above or below the scaling window. Above the window the
largest component is of order $\eps(n)n$ and it is the unique
component of this size. Below the window, the largest component is
of order $\eps^{-2}(n)\log(n\eps^3)$, but so is the $\ell$-th
largest component, for any fixed $\ell>1$. The following theorem
provides the analogous statement to (\ref{gnpsubcrit}) for
$G(n,d,p)$.
\begin{thm} \label{lower} {\bf [Below the critical window]}
Recall the definition of $\psi_n$ from (\ref{psidef}) and let
$\eps(n)>0$ be a sequence such that $\eps(n) \to 0$ and $\eps(n)
n^{1/3} \to \infty$. Consider $G(n,d,p)$ with $p = {1-\eps(n)
\over d-1}$ where $d\geq 3$ is fixed, then for any fixed integer
$\ell
> 0$ we have \be \label{subcrit} { |\C_\ell| \over
\psi_n(\eps(n))} \stackrel{P} \longrightarrow {d-2 \over d-1}
\qquad \hbox{as } n \to \infty \, , \ee
\end{thm}

\medskip

We now turn to the supercritical case. In the Erd\H{o}s-R\'enyi
random graph $G(n, {1 + \eps(n) \over n})$, where $\eps(n)>0$
satisfies $\eps(n) \to 0$ and $\eps(n) n^{1/3} \to \infty$ we have
(see \cite{B1}, \cite{L} and also \cite{NP2})
$$ { |\C_1| \over 2n\eps(n)} \stackrel{P} \longrightarrow 1 \qquad
\hbox{as } n \to \infty \, ,$$ and (\ref{gnpsubcrit}) holds for
any fixed integer $\ell > 1$, controlling the size of the smaller
components. The following theorem provides the analogous statement
for $G(n,d,p)$.

\begin{thm} \label{upper} {\bf [Above the critical window]}
Let $\eps(n)>0$ be a sequence such that $\eps(n) \to 0$ and
$\eps(n) n^{1/3} \to \infty$. Consider $G(n,d,p)$ with $p =
{1+\eps(n) \over d-1}$, where $d\geq 3$ is fixed, then
$$ { |\C_1| \over 2n\eps(n)} \stackrel{P} \longrightarrow {d \over d-2} \qquad
\hbox{as } n \to \infty \, .$$ Furthermore, for any fixed integer
$\ell > 1$ we have that (\ref{subcrit}) holds, controlling the
size of $C_\ell$.


\end{thm}

Next we turn to describe the limiting distribution of the
component sizes inside the scaling window $p={1+\lambda n^{-1/3}
\over d-1}$, in an analogous way to \cite{A}. Let $\{B(s) : s\in
[0,\infty)\}$ be standard Brownian motion and for $\lambda \in
{\mathbb R}$ define the process
\begin{equation}\label{process1} B^\lambda(s) = B\Big ( {(d-2)\over (d-1)}s \Big ) + \lambda s
- {(d-2) \over 2d}s^2 \, , \quad s\in [0,\infty) \, .
\end{equation}
Also, consider the reflected process
\begin{equation} \label{process2} W^\lambda(s) = B^\lambda(s) -
\min _{0 \leq s' \leq s} B^\lambda (s') \, .
\end{equation}
An {\rm excursion} $\gamma$ of $W^\lambda$ is a time interval
$[l(\gamma),r(\gamma)]$ in which $W^\lambda(l(\gamma)) =
W^\lambda(r(\gamma))=0$, and $W^\lambda(s) > 0$ for all $l(\gamma)
< s < r(\gamma)$. The excursion has length
$|\gamma|=r(\gamma)-l(\gamma)$. The sequence $(|\gamma _j|)_{j\geq
1}$ of excursion lengths, in decreasing order, is a random
variable in $\ell^2$ almost surely (see \cite{A}).

\begin{thm} \label{criticaldist} {\bf [Limiting distribution]}
Fix $\lambda \in {\mathbb R}$ and let $p={1+\lambda n^{-1/3} \over
d-1}$, where $d\geq 3$ is fixed, then
$$ n^{-2/3} \cdot(|{\cal C}_1|,|{\cal C}_2|, \ldots) \stackrel{d}{\Longrightarrow} (|\gamma _j|)_{j\geq
1}  \, ,$$ where convergence holds with respect to the $\ell^2$
norm.
\end{thm}

In \cite{NP3}, the authors prove that in bond percolation on any
$d$-regular graph on $n$ vertices with $p\leq {1 + \lambda
n^{-1/3}\over d-1}$, if the resulting graph typically has
components of size $n^{2/3}$ then their diameter is of order
$n^{1/3}$ and the mixing time of the lazy simple random walk on
these components is of order $n$. See \cite{NP3} for more details
and definitions. The following is an immediate corollary of
Theorem \ref{criticaldist} above and Theorem $1.2$ of \cite{NP3}.

\begin{cor} Consider $G(n,d,p)$ with $p={1+\lambda n^{-1/3} \over d-1}$ for
some $\lambda \in {\mathbb R}$, where $d\geq 3$ is fixed. Denote
by $\diam(\C_\ell)$ the diameter of $\C_\ell$ and let
$\Tm(\C_\ell)$ be the mixing time of the lazy simple random walk
on $\C_\ell$. Then for any fixed integer $\ell
> 0$ and any $\eps>0$ there exists $A=A(\eps, \lambda, \ell) <
\infty$ such that for all large $n$,
\begin{itemize}
\item $ \quad \prob \Big ( \diam(\C_\ell)\not \in [A^{-1} n^{1/3},
An^{1/3}] \Big ) < \eps \, ,$ \item $ \quad \prob \Big (
\Tm(\C_\ell) \not \in [A^{-1} n, An] \Big ) < \eps \, .$
\end{itemize}
\end{cor}


A major challenge is to give criteria for specific $d$-regular
graphs to exhibit mean-field behavior (the theorems of this paper
establish that this occurs for most $d$-regular graphs).
Substantial progress in this direction was made in \cite{BCHSS1}
and
\cite{BCHSS2}.\\

The rest of the paper is organized as follows. As the proof of
Proposition \ref{dregular} is simple and instructive, we provide
it in Section $2$. In Section $3$ we describe a discrete
exploration process which generates a random sample of $G(n,d,p)$.
The analysis of this process is crucial for proving the results of
this paper and is presented in Section $4$. From there we proceed
to prove Theorem \ref{criticalupper} in Section $5$. Theorems
\ref{lower} and \ref{upper}, describing the behavior above and
below the scaling window, are proved in Section $6$ and $7$
respectively. Theorem \ref{criticaldist} is proved in Section
\ref{seclimit} and we end with some concluding remarks in Section
\ref{secconc}. \\

We use the standard asymptotic notation. For two functions $f(n)$
and $g(n)$, we write $f=o(g)$ if $\lim_{n\rightarrow\infty}
f/g=0$. Also, $f=O(g)$ if there exists an absolute constant $C>0$
such that $f(n)<Cg(n)$ for all large enough $n$ and $f=\Theta(g)$
if both $f=O(g)$ and $g=O(f)$ hold.

\section{Proof of the general upper bounds (Prop. \ref{dregular})}

For the proofs in this section and in sections to follow we
present some standard facts about processes with independent
increments.

\begin{lemma} \label{gambler} Let $\beta$ be a random variable
supported on the integers with $\prob ( \beta < -1) = 0$. Let
$\{\beta_i\}$ be i.i.d. random variables distributed as $\beta$
and let $W_t = W_0 + \sum_{i=1}^t \beta_i$, where $W_0>0$ is some
integer. For an integer $h>W_0$ define the stopping time
$$ \gamma_h = \min _t \{ W_t = 0 \textrm{ or } W_t \geq h \} \,
.$$ We have
\begin{enumerate}
\item[(i)] If $c>0$ is such that $\E e^{-c\beta} \geq 1$ then
$$ \prob (W_{\gamma_h} > 0) \leq {1- e^{-cW_0} \over 1-e^{-ch}} \,
.$$

\item[(ii)] If $c>0$ is such that $\E e^{c\beta} \leq 1$ then
$$ \prob (W_{\gamma_h} > 0) \leq {e^{cW_0} - 1 \over e^{ch} -1 }
\, .$$
\end{enumerate}
\end{lemma}

\noindent {\bf Proof.} This is a standard application of the
optional stopping theorem (see \cite{D}). The assumption $\E
e^{-c\beta} \geq 1$ implies that $\{e^{-cW_t} \}$ is a
submartingale. Optional stopping gives
$$ e^{-cW_0} \leq 1 - \prob(W_{\gamma_h} > 0) + e^{-ch}\prob(W_{\gamma_h} >
0) \, ,$$ which yields assertion (i) of the lemma. The assumption
$\E e^{c\beta} \leq 1$ implies that $\{ e^{cW_t} \}$ is a
supermartingale, and similarly we get assertion (ii) of
the lemma. \qed \\

The following lemma is a variant of a lemma due to Bahadur and
Rao~\cite{BR}.

\begin{lemma} \label{tail} Let $\beta$ be a non-lattice, integer valued
random variable with $\E \beta^2 < \infty$. Let $\{\beta_i\}$ be
i.i.d. random variables distributed as $\beta$ and let $W_t = W_0
+ \sum_{i=1}^t \beta_i$. Let $\tau$ be the hitting time of $0$,
i.e.
$$ \tau = \min _t \{ W_t = 0\} \, .$$
If $\theta_0>0$ satisfies

\be \label{tailassump1} \E [\beta e^{\theta_0 \beta}] = 0 \, ,\ee
then for any integer $\ell>0$ we have
$$ \prob ( \tau = \ell) = \Theta \Big (\ell^{-3/2} \varphi(\theta_0) ^\ell \Big
)\, ,$$
where $\varphi(\theta) = \E e^{\theta \beta}$, and the constants
in the $\Theta$ depend only on $\beta$ and $W_0$ but not on
$\ell$.
\end{lemma}

For the proof of Lemma \ref{tail} we require the following variant
of a lemma due to Spitzer \cite{S}. For completeness, we include
its proof here.
\begin{lemma} \label{spitzer}
Let $a_0, \ldots, a_{k-1} \in {\mathbb Z}$ be such that $\sum
_{i=0}^{k-1} a_i = -d$. Then there are at least one and at most
$d$ numbers $j \in \{0, \ldots, k-1\}$ such that for all $\ell \in
\{0, \ldots, k-2\}$
$$ \sum _{i=0}^{\ell} a_{(j+i) {\rm \  mod \ } k} > -d \, .$$
\end{lemma}
\noindent {\bf Proof.} Continue the sequence periodically such
that $a_{k+s} = a_{s}$ for any integer $s > 0$. Let $j$ be the
first global minimum of the function $f(j) = \sum_{i=0}^j a_i$ on
the domain $\{0,\ldots, k-1\}$. It is easy to see that for that
$j$, and any $\ell \in \{0, \ldots, k-2\}$
$$ \sum _{i=0}^{\ell} a_{(j+i)} > -d \, .$$
Assume now that there were $j_1 < \ldots < j_{d+1}$ all in
$\{0,\ldots ,k-1\}$ satisfying that for all $\ell \in \{0, \ldots,
k-2\}$
$$ \sum _{i=0}^{\ell} a_{(j_r + i)} > -d \, ,$$
for all $r \in \{1, \ldots, d+1\}$. Define a function $g(r)$ on
$\{1, \ldots, d+1\}$ by
$$ g(1) = \sum _{i=j_{d+1}}^{k-1+j_1 -1} a_i\, ,$$
$$ g(r) = \sum _{i=j_{r-1}}^{j_r -1} a_i \, , \quad r \in \{2,\ldots,d+1\} \, .$$
As $\sum_ {i=j_{r-1}}^{k-1+j_{r-1}} a_i = -d$ and $\sum_ {i=j_r
}^{k-1+j_{r-1}} a_i > -d$ we find that $g(r) \leq -1$ for all $r
\in \{1, \ldots, d+1\}$. The assumption $\sum _{i=0}^{k-1} a_i =
-d$ implies that $g(1)+\ldots+g(d+1)=-d$ and we have arrived
at a contradiction. \qed \\

\noindent {\bf Proof of Lemma \ref{tail}.} Let $\beta_\theta$ be a
random variable distributed as
$$ \prob (\beta_\theta = t) = \varphi(\theta)^{-1} e^{\theta t} \prob
(\beta = t) \, .$$ Let $\{ \beta_\theta(i)\}$ be a sequence of
i.i.d. random variables distributed as $\beta_\theta$ and let
$W_\theta(t) = W_0 + \sum_{i=1}^t \beta_\theta(i)$.
Let $I = \{ (t_1, \ldots, t_\ell) : W_0 + t_1 + \ldots + t_\ell =
0 \}$ and observe that
$$ { \prob (W_\ell = 0) \over \prob (W_\theta(\ell) = 0)} = { \sum
_{(t_1, \ldots, t_\ell) \in I} \prod _{i=1}^\ell \prob ( \beta_i =
t_i) \over \varphi(\theta)^{-\ell} \sum _{(t_1, \ldots, t_\ell)
\in I} e^{\theta(t_1 + \ldots + t_\ell)} \prod _{i=1}^\ell \prob (
\beta_i = t_i) } = e^{\theta W_0} \varphi(\theta)^\ell .$$ We now
take $\theta = \theta_0$. By (\ref{tailassump1}) we have that $\E
\beta_{\theta_0} = 0$, thus by the local central limit theorem
(see \cite{D}, Section 2.5) we have that $\prob
(W_{\theta_0}(\ell) = 0) = \Theta( \ell ^{-1/2} )$. Thus,
$$ \prob (W_\ell = 0) = \Theta (\ell ^{-1/2} \varphi(\theta_0)^\ell)
\, ,$$ and by Lemma \ref{spitzer} we learn that $\prob (\tau =
\ell) = \Theta (\ell^{-1}) \prob (W_\ell = 0)$, concluding
our proof. \qed \\

\noindent {\bf Proof of Proposition \ref{dregular}.} For a graph
$G$, denote by $G_p$ the random graph obtained by bond percolation
on $G$ with probability $p$. For a vertex $v$ and let $\C(v)$
denote the connected component that contains $v$ in $G_p$. We
recall an exploration process, developed independently by
Martin-L\"of \cite{M} and Karp \cite{K}. In this process, vertices
will be either {\em active, explored} or {\em neutral}. At each
time $t$, the number of active vertices will be denoted $Y_t$ and
the number of explored vertices will be $t$. Fix an ordering of
the vertices, with $v$ first. As an upper bound, assume some edge
$(v,u)$ adjacent to $v$ is open. At time $t=0$, the vertices $v$
and $u$ are active and all other vertices are neutral, so $Y_0=2$.
In step $t>0$ let $w_t$ be the first active vertex. Denote by
$\eta_t$ the number of neutral neighbors of $w_t$ in $G_p$ and
change the status of these vertices to {\em active}. Then, set
$w_t$ itself {\em explored}. The process stops when $Y_t$ hits
$0$, and observe that since at each step we set precisely one
vertex explored we have $|\C(v)|\leq \min \{t : Y_t=0\}$. Let
$\{w_1, w_2, \ldots\}$ be independent random variables distributed
as Bin$(d-1,p)-1$. Let $W_t = 2+\sum _{i=1}^t w_i$. As $G$ is
$d$-regular, it is clear that we can couple the process $\{Y_t\}$
and $\{W_t\}$ such that $Y_t \leq W_t$ for all $t\leq
|\C(v)|$. \\

We begin with the proof of part $1$ of the Theorem. We will use
Lemma \ref{tail} with $\beta = w - 1$, where $w$ is distributed as
Bin$(d-1,p)$  and $p={1- \eps \over d-1}$. If we write $w = \sum
_{j=1}^{d-1} I_j$ where $I_j$ are i.i.d. Bernoulli$(p)$ random
variables, we get that for any $\theta>0$
$$ \E we^{\theta w} = (d-1)\E \Big [ \prod _{j=2}^{d-1} e^{\theta
I_j} \Big ] \E I_1 e^{\theta I_1} = (d-1)pe^{\theta}
(1-p+e^{\theta}p)^{d-2} \, .$$ As $\E e^{\theta w} =
(1-p+pe^\theta)^{d-1}$ we have
$$ \E \beta e^{\theta \beta} = e^{-\theta}
(1-p+e^{\theta}p)^{d-2}  [ p(d-1)e^{\theta} - (1-p+pe^{\theta}) ]
\, .$$ Let $\theta_0>0$ be a number such that $\E \beta
e^{\theta_0 \beta} = 0$, then by estimating $e^x = 1+x+O(x^2)$ in
the last equation we find that
$$ p(d-2)(1+\theta_0) + p + O(\theta_0 ^2) = 1 \, ,$$
thus
$$ \theta_0 = { (d-1)\eps \over d-2} + O(\eps^2) \, .$$
For any $\theta>0$ by estimating $e^x = 1+x+x^2/2 + O(x^3)$ we get
\begin{eqnarray*} \varphi(\theta) &=&
\E e^{\theta \beta} = e^{-\theta} ( 1 + p(e^{\theta} -1)) ^{d-1}
 \\ &=& \Big (1-\theta + \theta^2/2 \Big )\Big (1+ (d-1)p(\theta + \theta^2/2) + {(d-1)(d-2) p^2 \theta^2 \over 2}\Big ) + O(\theta^3)  \, .
\end{eqnarray*}
By simplifying and plugging in the value of $\theta_0$ we find
that
$$ \varphi(\theta_0) = 1 - {(d-1) \eps ^2  \over 2(d-2)} + O(\eps^3) \, .$$
Let $\tau = \min\{t : W_t = 0\}$ then Lemma \ref{tail} implies
that
$$ \prob ( \tau > T) \leq \sum _{\ell =T+1}^\infty O \Big (
\ell ^{-3/2} \Big (1 - {(d-1) \eps ^2  \over 2(d-2)} + O(\eps^3)
\Big )^\ell \Big ) \, .$$ We take
$$ T = (1+\eta) {2(d-2)\over d-1} \eps^{-2} \log(n\eps^3) \, ,$$
and a straightforward computation using $1-x \leq e^{-x}$ yields
that for some fixed $c>0$
$$ \prob (\tau > T) \leq O\Big ( \eps
(n\eps^3)^{-(1+\eta)(1-c\eps)} \log (n\eps^3)^{-3/2} \Big ) \, .$$
As $Y_t \leq W_t$ for all $t\leq |\C(v)|$ we have $\prob( |\C(v)|
> T) \leq \prob(\tau
> T)$. Denote by $X$ the number of vertices $v$ of $G$ such that
$|\C(v)|>T$. If $|\C_1| > T$ then $X>T$. We conclude that for some
$c_1>0$

\begin{eqnarray*} \label{useful} \prob ( |\C _1 | > T) &\leq& \prob (X > T) \leq
{\E X \over T} \leq { n \prob (|C(v_1)| > T) \over T} \\ &\leq& {C
n \eps (n\eps^3)^{-(1+\eta)(1-c\eps)} \over \eps^{-2}
\log^{5/2}(n\eps^3)} \leq (n \eps ^3)^{{-\eta}(1-c_1
\eps)+c_1\eps} \to 0  \, ,
\end{eqnarray*}
which concludes part $1$ of the proposition. \\

We now prove part $2$ of the Proposition, following the strategy
laid out in \cite{NP}. By monotonicity we may assume that $p={1
\over d-1}$. In that case $\{W_t\}$ is a martingale with $\E W_0
\leq 2$. Define $\gamma_h$ as in Lemma \ref{gambler}, so by
optional stopping we get that $\E W_0 = \E W_{\gamma_h} \geq h
\prob (W_{\gamma_h}
>0)$, whence
\be \label{part1equ1} \prob (W_{\gamma_h} >0) \leq {2 \over h} \,
.\ee By Corollary $6$ in \cite{NP} (see also inequality ($3$) of
\cite{NP}) we also have

\be \label{part1equ2} \E [ W_{\gamma_h}^2 \mid W_{\gamma_h} >0 ]
\leq h^2 + 3h \, .\ee It is immediate to verify that $W_t^2 - (1-
{1 \over d-1})t$ is also a martingale. Optional stopping,
(\ref{part1equ1}) and (\ref{part1equ2}) gives that
$$  (1- {1 \over d-1})\E \gamma_h \leq \E W_{\gamma_h}^2 = \prob
(W_{\gamma_h} >0) \E [ W_{\gamma_h}^2 \mid W_{\gamma_h} >0 ] \leq
2h+6 \, .$$ As $d>2$ we get

\be \label{part1equ3} \E \gamma_h \leq 4h + 12 \, .\ee Hence as
long as $h>12$
$$ \prob (\gamma_h \geq h^2) \leq {5 \over h} \, .$$
Denote $\gamma_h^* = \gamma_h \wedge h^2$. By the previous
inequality and (\ref{part1equ1}), we have
$$ \prob ( W_{\gamma_h^*} > 0 ) \leq \prob ( W_{\gamma_h} > 0 ) +
\prob ( \gamma_h \geq h^2 ) \leq {7 \over h} \, .$$ Let $T=h^2$
and observe that if $|\C(v)| > h^2$ we must have $W_{\gamma_h^*}
>0$, thus $\prob ( |\C(v)|>T ) \leq {7 \over \sqrt{T}}$. Again denote
by $X$ the number of vertices $v$ of $G$ such that $|\C(v)|>T$. If
$|\C_1| > T$ then $X>T$. We put $T= \Big ( \lfloor \sqrt{An^{2/3}}
\rfloor \Big )^2$ and conclude that
$$ \prob ( |\C_1| > T ) \leq {\E X \over T} \leq { 7 n \over
T^{3/2} } \leq {8 \over A^{3/2}} \, ,$$ for large enough $n$, as required.\\

We now prove part $3$ of the Theorem. For an integer $k>0$ denote
by $X_k$ the number of vertices $v$ of $G$ such that $|\C(v)|> k$.
It is clear that if $|\C_1| > k$ then $|\C_1| \leq X_k$, thus \be
\label{part3} \E |\C_1(G_p)| \leq k + \E X_k \, .\ee We estimate
the last term of the previous display in a similar way to the
proof of part $1$ of the proposition. Put $p={1 + \eps \over
d-1}$, let $w$ be distributed as Bin$(d-1,p)$ and $\beta = w -1$.
By an almost identical calculation to the one done in part $1$ we
get that in the notation of Lemma \ref{tail}
$$ \theta_0 = - {(d-1)\eps \over d-2} + O(\eps^2) \, ,$$
and
$$ \varphi(\theta_0) = 1 - {(d-1) \eps ^2  \over 2(d-2)} + O(\eps^3) \, .$$
Lemma \ref{tail} and our usual coupling gives that for some $C>0$,
$$ \prob ( |\C(v)| > k ) \leq \prob ( \tau > k ) \leq \sum _{\ell
> k} C \ell^{-3/2} (1- \Theta(\eps^2)) ^\ell \, .$$
A straightforward calculation with the sum in the previous display
shows we can bound it from above by $C(k^{-1/2} + \eps)$ for some
fixed $C>0$. We find that $\E X_k = n \prob ( |\C(v)| > k ) \leq C
\eps n + C n k^{-1/2}$. Choosing $k=n^{2/3}$ and plugging into
(\ref{part3}) concludes the proof. \qed \\

\section{The random regular graph and the exploration process} \label{explore}

The following model, known as the {\em configuration model}, was
introduced by Bollob\'as in \cite{B0} (see also \cite{BC} and
\cite{Wor}) and was used to construct a uniform random $d$-regular
graph on $n$ vertices, assuming $dn$ is even. Consider the vertex
set $\{1, \ldots , dn\}$ as $n$ distinct $d$-tuples. Draw a
uniform perfect matching on the set $\{1, \ldots , dn\}$, and then
contract every $d$-tuple into a single vertex. It was shown in
\cite{B0} and \cite{BC} that with probability tending to
$\exp({1-d^2 \over 4})$ as $n\to \infty$ this process yields a
simple $d$-regular graph. Moreover, conditioning on this event,
the graph obtained is uniformly distributed among all {\em simple}
$d$-regular graphs on $n$ vertices.

A uniform perfect matching on a set can be obtained by drawing the
edges of the matching sequentially: for each edge choose the first
vertex according to any rule (deterministic or random) and then
choose the second vertex uniformly at random among the unmatched
vertices.
This motivates exploring the connected components (in the spirit
of \cite{K} and \cite{M}) by drawing a uniform matching on
$\{1,\ldots ,dn\}$ sequentially, and independently percolating
each edge of the matching; we call this process the {\em
exploration process}. In this process, vertices will be either
{\em active, explored} or {\em neutral} and each $d$-tuple may
contain vertices with different status. Choose an ordering of the
vertices $\{ v_{i,k} \, : \, 1 \leq i \leq n \, , 1 \leq k \leq d
\}$ where $\{v_{i,1}, \ldots , v_{i,d} \}$ is the $i$-th
$d$-tuple, for $0 \leq k \leq n-1$. Initially, the first
$d$-tuple, vertices $\{v_{1,1}, \ldots, v_{1,d}\}$, are active and
all other vertices are neutral. At each time $t > 0$, if there are
active vertices, let $w_t$ be the first active vertex; if there
are no active vertices, let $w_t$ be the next neutral vertex and
change the status of the neutral vertices in $w_t$'s $d$-tuple to
{\em active} (including the status of $w_t$ itself). Now match
$w_t$ with a uniformly drawn unmatched vertex $\eta_t$.
If $\eta_t$ is neutral and the edge $(w_t, \eta_t)$ is retained in
the percolation then we change the status of the neutral vertices
in $\eta_t$'s $d$-tuple to {\em active}, and we also set $w_t$ and
$\eta_t$ {\em explored}. If $\eta_t$ is neutral and the edge
$(w_t, \eta_t)$ is not retained in the percolation or if $\eta_t$
is active, just set $w_t$ and $\eta_t$ {\em explored} without
changing the status of any other vertex. This gives a graph on
$\{v_1, \ldots, v_{dn}\}$; we obtain the multi-graph $G^*(n,d,p)$
on $n$ vertices by contracting each $d$-tuple to a single vertex.
Denote by $\Simp$ the event that the perfect matching constructed
by the exploration process yields a simple $d$-regular graph. By
\cite{B0} and our previous discussion we have

\be \label{simple} \prob ( \Simp ) = \exp \Big ({1-d^2 \over 4}
\Big )+o(1) \, , \ee and by our previous discussion, if we
condition on this
event, then $G^*(n,d,p)$ is distributed as $G(n,d,p)$. \\

In order to analyze the exploration process we introduce the
following random variables. For $0 \leq k \leq d$ and $t \leq
dn/2$ denote by $\Ne^{(k)}_t$ the set of $d$-tuples which have
precisely $k$ neutral vertices {\em after} $\eta_t$ was drawn and
{\em before} $w_{t+1}$ is chosen, and by $\Net^{(k)}_t$ the set of
$d$-tuples which have precisely $k$ neutral vertices {\em after}
$w_{t+1}$ was chosen and {\em before} $\eta_{t+1}$ is drawn. Let
$N^{(k)}_t$ and $\NN^{(k)}_t$ denote the cardinality of these
sets, respectively. For a vertex $v$, denote by $[v]$ the tuple
containing $v$. Hence, the notation $[w_{t+1}] \in \Ne^{(k)}_{t}$
implies that the $d$-tuple of $w_{t+1}$, after $w_{t+1}$ was
chosen, has precisely $k$ neutral vertices, and therefore
$w_{t+1}$ was chosen neutral (i.e., there were no active vertices
remaining). Similarly, $[\eta_t] \in \Net^{(k)}_{t-1}$ is the {\em
event} that the $d$-tuple of $\eta_t$ has $k$ neutral vertices
after $\eta_t$ was drawn, and that $\eta_t$ was drawn neutral. For
an edge $e$ we write $e \in G_p$ to denote that $e$ was retained
in the percolation.

The exploration process dictates the recursive dynamics of these
random variables. The number of $d$-tuples which have $d$ neutral
vertices after $w_1$ is chosen is $n-1$; at each time $t>0$ we
have $N^{(d)}_t = \NN^{(d)}_{t-1} - 1$ if $[\eta_t] \in
\Net^{(d)}_{t-1}$ and $\NN^{(d)}_t = N^{(d)}_{t} -1$ if $[w_{t+1}]
\in \Ne^{(d)}_{t}$. Hence,
$$ \NN^{(d)}_0 = n-1 \, ,$$

\be \label{dneutral} N^{(d)}_t = \NN^{(d)}_{t-1} - {\bf 1}_{\{
[\eta_t] \in \Net^{(d)} _{t-1}\}} \, , \quad \NN^{(d)}_t =
N^{(d)}_{t} - {\bf 1}_{\{ [w_{t+1}] \in \Ne^{(d)} _{t}\}}  \, .\ee

For $0 < k < d$, at time $t=0$ there are no $d$-tuples with $k$
neutral vertices. At each time $t>0$ we have $N^{(k)}_t =
\NN^{(k)}_{t-1} - 1$ if $[\eta_t] \in \Net^{(k)}_{t-1}$ and
$N^{(k)}_t = \NN^{(k)}_{t-1} + 1$ if $[\eta_t] \in
\Net^{(k+1)}_{t-1}$ and the edge $(w_t,\eta_t)$ is {\em not}
retained in the percolation. We also have $\NN^{(k)}_t =
N^{(k)}_{t} - 1$ if $[w_{t+1}] \in \Ne^{(k)}_{t}$. Hence,

$$ \NN^{(k)}_0  = 0 \, ,$$

\be \label{kneutral} N^{(k)}_t = \NN^{(k)}_{t-1} - {\bf 1}_ {\{
[\eta_{t}] \in \Net^{(k)} _{t-1}\}} +  {\bf 1} _{\{(w_{t},
\eta_{t}) \not \in G_p\}} {\bf 1}_ {\{ [\eta_{t}] \in \Net^{(k+1)}
_{t-1}\}} \, . \ee

\be \label{kneutral2} \NN^{(k)}_t =  N^{(k)}_{t} - {\bf 1}_ {\{
[w_{t+1}] \in \Ne^{(k)} _{t}\}} \, , \ee

Finally we have $\NN^{(0)}_0=1$ (as the $d$-tuple of $w_1$ has no
neutral vertices) and at each time $t>0$ we have $N^{(0)}_t =
\NN^{(0)}_{t-1} + 1$ if $\eta_t$ is drawn neutral and the edge
$(w_{t}, \eta_{t})$ was retained in the percolation, or if
$[\eta_{t}] \in \Net^{(1)}_{t-1}$ and the edge $(w_{t}, \eta_{t})$
was {\em not} retained in the percolation. We also have
$\NN^{(0)}_t = N^{(0)}_{t} + 1$ if $w_{t+1}$ is chosen neutral,
i.e., in the case where no more active vertices are left. Hence,
$$\NN^{(0)}_0=1 \, , $$
\begin{eqnarray}
 \label{0neutral} N^{(0)}_t = \NN^{(0)}_{t-1} &+& {\mathbf 1} _{\{(w_{t}, \eta_{t}) \in G_p\}} {\mathbf
1}_{\{\eta_{t} \hbox{ drawn neutral}\}} \nonumber \\ &+&  {\mathbf
1} _{\{(w_{t}, \eta_{t}) \not \in G_p\}} {\mathbf 1}_{\{[\eta_{t}]
\in \Net_{t-1}^{(1)} \}} \, ,
\end{eqnarray}

\be \label{0neutral2} \NN^{(0)}_t = N^{(0)}_{t} + {\mathbf 1}_{ \{
w_{t+1} \hbox{ chosen neutral} \}} \, . \ee

Denote by $\Ac_t$ the set of active vertices {\em after} $\eta_t$
was drawn and {\em before} $w_{t+1}$ is chosen and by $\Act_t$ the
set of active vertices {\em after} $w_{t+1}$ was chosen and {\em
before} $\eta_{t+1}$ is drawn. Let $A_t$ and $\AA_t$ denote the
cardinality of these sets, respectively. Let $\{\xi_t\}$ be random
variables defined by

\be \label{incdef} \xi_t = {\mathbf 1} _{\{(w_t, \eta_t) \in
G_p\}} \sum_{k=2}^d (k-1) {\mathbf 1}_{\{[\eta_t] \in \Net^{(k)}
_{t-1}\}} - {\mathbf 1}_{\{ \eta_t \in \Act_{t-1}\}} -1 \, . \ee
For the vertex $w_t$ denote by $N(w_t)$ the number of neutral
vertices in $[w_t]$ {\em after} $w_t$ was chosen and {\em before}
$\eta_t$ was drawn, including $w_t$ itself. Note that if $w_t$ is
active then $N(w_t) = 0$, so this number is non-zero only if $w_t$
is neutral, i.e., when $A_{t-1}=0$.

We now describe the recursive dynamics of these random variables.
After choosing $w_1$ and before choosing $\eta_1$ we have
precisely $d$ active vertices hence $\AA_0=d$. If $[\eta_t] \in
\Net^{(k)} _{t-1}$ and the edge $(w_t,\eta_t)$ was retained in the
percolation then we mark $k-1$ neutral vertices as active
vertices, and one active vertex as explored, so $A_t = \AA_{t-1} +
(k-1) - 1$. Also, if $[\eta_t] \in \Net^{(k)} _{t-1}$ but the edge
$(w_t,\eta_t)$ was {\em not} retained in the percolation then $A_t
= \AA_{t-1} - 1$. If $\eta_t \in \Act_{t-1}$ then we mark two
active vertices as explored and hence $A_t = \AA_{t-1} - 2$.
Together this gives
$$ A_t = \AA_{t-1} + {\mathbf 1} _{\{(w_t, \eta_t) \in
G_p\}} \sum_{k=2}^d (k-1) {\mathbf 1}_{\{[\eta_t] \in \Net^{(k)}
_{t-1}\}} - {\mathbf 1}_{\{ \eta_t \in \Act_{t-1}\}} -1 \, .$$

If $A_t > 0$ then $w_{t+1}$ will be chosen active and so $\AA_t =
A_t$. On the other hand, if $A_t = 0$ then we mark the neutral
vertices in $[w_{t+1}]$ (including $w_{t+1}$ itself) as active,
and hence $\AA_t = N(w_{t+1})$. This together with the previous
display and (\ref{incdef}) gives

$$ A_0 = 0 \, , \qquad \AA_0 = d \, ,$$

\be \label{arec} A_t = A_{t-1} + \xi_t + N(w_t)  \, . \ee \\


Let $0 = t_0 < t_1 < t_2 < \ldots$ be the times at which
$A_{t_j}=0$. At time $t_j$ we completely explored the $j$-th
component (which we have started exploring in time $t_{j-1}+1$)
and all the $d$-tuples that became completely explored between
times $t_{j-1}+1$ and $t_j$ are the vertices of this component.
Define the random variables
$$ S_j = \Big | \Big \{ t \in (t_{j-1}, t_j] : (w_t, \eta_t) \in
G_p \hbox{ and } \eta_t \hbox{ is drawn neutral} \Big \} \Big | \,
,$$
$$ U_j = \Big | \Big \{ t \in (t_{j-1}, t_j] : \eta_t \in \Act_{t-1}\Big \} \Big | \, ,$$
$$ V_j = \Big | \{ t \in (t_{j-1}, t_j] : [\eta_t] \in
\cup_{i=1}^{d-1} \Net_{t}^{(i)} \Big \} \Big | \, .$$

The following lemma relates all the above to component sizes of
the graph $G^*(n,d,p)$.

\begin{lemma}\label{bottomline} The size of the $j$-th completely explored component is
$S_j+1$. Furthermore, we have
$$ 0 \leq S_j + 1 - {t_j - t_{j-1} \over d-1}
\leq {U_j \over d-1} + V_j+1 \, . $$
\end{lemma}
\noindent {\bf Proof.} At each time where $\eta_t$ is neutral and
$(w_t, \eta_t)\in G_p$ we add a new $d$-tuple to our currently
explored component, increasing its size by $1$. Thus, the size of
the $j$-th completely explored component is simply $S_j+1$. To get
the second part of the lemma denote by $T_j$ the random variable
$$ T_j = \Big | \Big \{ t \in (t_{j-1}, t_j] : (w_t, \eta_t) \not \in
G_p \hbox{ and } \eta_t \hbox{ is drawn neutral} \Big \} \Big | \,
.$$ Observe that since $\eta_t$ is drawn among the neutral and
active vertices remaining we have \be \label{eq1} t_{j} - t_{j-1}
= S_j + T_j + U_j \, .\ee Consider now the dynamics described in
the two paragraphs preceding (\ref{arec}). By the previous
display, and since $A_{t_{j-1}}=A_{t_{j}} = 0$ we have
$$ 0 = N(w_{t_{j-1}+1})-1 -2 U_j - T_j +
\sum _{k=1}^{d} (k-2) \cdot \Big | \Big \{ t \in (t_{j-1},t_j] :
(w_t, \eta_t) \in G_p \hbox{ and } [\eta_t] \in \Net_{t}^{(k)}
\Big \} \Big | \, .$$ The last sum in the equation can be bounded
above by $(d-2)S_j$ and below by $(d-2)S_j - (d-1)V_j$. This
together with (\ref{eq1}) and the fact that $1 \leq
N(w_{t_{j-1}+1}) \leq d$ gives that
$$ 0 \leq S_j + 1 - {t_j - t_{j-1} \over d-1}
\leq {U_j \over d-1} + V_j+1 \, .$$ \qed \\

It will be more convenient to work with the process $\{Y_t\}$
defined by
$$ Y_0 = d\, , \qquad Y_t = Y_{t-1} + \xi _t \, .$$ There is an evident
connection between the process $\{Y_t\}$ and $\{A_t\}$. By
(\ref{arec}) we have
$$ Y_t = A_t - Z_t \, ,$$
where
$$Z_t = \sum _{i=1}^t N(w_i) \, .$$
Observe that $Z_t$ is an increasing process and $Z_t=Z_{t_j+1}$
for all $t \in \{t_j+1, \ldots, t_{j+1}\}$. As $A_{t_j}=0$ we have
that $Y_{t_j} = -Z_{t_j}$ for all $j$. Thus, for any $t \in \{t_j
+1, \ldots, t_{j+1}-1\}$ we have

$$ Y_{t_{j+1}} = -Z_{t_{j+1}} = -Z_t  < Y_t \, ,$$
as $A_t > 0$ for such $t$'s. By induction we learn that
$Y_{t_{j+1}} < Y_t$ for all $t < t_{j+1}$. Hence, the $t_j$'s are
 record minima for the process $\{Y_t\}$. Since $N(w_{t_j+1}) \leq d$
 we have that $Z_{t_j+1} \leq -Y_{t_j+1} + d$. Thus, by our
 previous discussion we learn that $Z_t \leq -\min_{s \leq t}Y_s
 +d$. We conclude that
\be \label{activecount2} A_t \leq Y_t - \min \{ Y_s : s \leq t\} +
d \, .\ee

\section{Exploration Process Analysis}

For the following, we assume that $\eps=\eps(n)$ is a sequence
such that $\eps(n) \to 0$ and we write $p=p(n)={1+\eps(n)\over
d-1}$. Let ${\cal F}_t$ be the $\sigma$-algebra
$$ \F_t = \sigma \Big \{ N^{(k)}_j\, , \, \NN^{(k)}_j :
0\leq j \leq t \, , \,  0\leq k \leq d \Big \} \, .$$
At each time $t$ we have that $\eta_t$ is chosen uniformly among
the $dn-2t+1$ neutral and active vertices remaining (which are not
$w_t$). Thus for any $0 \leq k \leq d$ we have

\be \label{uniform0} \E \Big [ {\mathbf 1}_{\{[\eta_t] \in
\Net^{(k)} _{t-1}\}} \mid {\cal F}_{t-1} \Big ] = {k
\NN_{t-1}^{(k)} \over dn-2t+1}\, , \ee and

\be \label{uniform2} \quad \E \Big [ {\mathbf 1}_{\{\eta_t \in
\Act _{t-1}\}} \mid {\cal F}_{t-1} \Big ] = {\AA_{t-1} \over
dn-2t+1} \, , \ee hence

\be \label{uniform} \prob \Big ( [\eta_t] \in \Net^{(k)} _{t-1}
\Big ) = {k \E \NN_{t-1}^{(k)} \over dn-2t+1}  \, , \quad \prob
\Big ( \eta_t \in \Act _{t-1} \Big ) = {\E \AA_{t-1} \over
dn-2t+1}  . \ee \\

%

In the conditions of Lemmas \ref{activeest} - \ref{deviation}
below and Corollary \ref{drift} appears a constant $C$ and the
constants implicit in the $O$-notation depend on $C$.


\begin{lemma} \label{activeest} For any $C>0$ we have that for all
$t<C\eps(n)n$

\be \label{aest} \E A_t = O(\eps t + \sqrt{t}) \, ,\ee and

\be \label{zest} \E Z_t = O(\eps t + \sqrt{t}) \, .\ee
\end{lemma}

\begin{lemma} \label{moments} For any $C>0$ we have that for all
$t<C\eps(n)n$

\be \label{dbound}  \E \NN_t^{(d)} = n-t +  O(\eps t+\sqrt{t}) \,
, \ee

\be \label{d-bound} \E \NN_{t}^{(d-1)} = (1-p)t + O(\eps
t+\sqrt{t}) \, , \ee


\be \label{kbound} \E \NN_{t}^{(k)} = O(\eps t) \, , \qquad  0 < k
< d-1 , \ee
\end{lemma}

\begin{lemma} \label{deviation} For any $C>0$ we have that for all
$t<C\eps(n)n$

\be \E | \NN_t^{(k)} - \E \NN_t^{(k)} | \leq O(\eps t + \sqrt{t})
\, , \qquad 0 \leq k \leq d \, .  \ee
\end{lemma}


\begin{cor}\label{drift} For any $C>0$ we have that for all
$t<C\eps(n)n$
\begin{enumerate}
\item[(i)] $  \E \xi_t  - \eps + {d-2 \over d(d-1)} \cdot {t
\over n} = O \Big ( \eps^2 + {\sqrt{t} \over n} \Big ) \, ,$
\item[(ii)] $ \E \Big | \E [ \xi_t \mid \F_{t-1} ] - \E \xi_t  \Big | = O\Big ( {\eps
t + \sqrt{t} \over n} \Big ) \, ,$
\item[(iii)] $ \E [ \xi_t^2 \mid \F_{t-1} ]  - (d-2) = O(\eps) \, .$

\end{enumerate}
\end{cor}


\begin{lemma} \label{extradev}
For any small $\delta>0$ there exists some constant
$c=c(\delta)>0$ such that if $t \leq \delta n$ then  \be
\label{deviation2} \prob \Big ( \NN_{t}^{(d)}
> n- (1-3\delta)t \Big ) \leq e^{-ct} \, ,\ee and \be
\label{deviation3} \prob \Big ( \NN_{t}^{(0)} < pt(1-3\delta) \Big
) \leq e^{-ct} \, . \ee
\end{lemma}

In order to bound the terms $U_j$ and $V_j$ in Lemma
\ref{bottomline} we have the following lemma.

\begin{lemma} \label{excessbd} For an integer $0<T< n/4$ define
$$ U_T = \Big \{ t \leq T : \eta_t \in \Act_{t-1} \hbox{ or } \eta_t \in \cup_{i=1}^{d-1} \Net_{t-1}^{(i)} \Big \} \, .$$ Then there exists some constant
$c>0$ such that if $4\sqrt{n}< T < n/4$
$$ \prob \Big ( |U_T| > {4 T^2 \over n} \Big ) \leq e^{-cT^2 /n} \, .$$
\end{lemma}


%

\noindent {\bf Proof of Lemma \ref{activeest}.} We rely on the
inequality (\ref{activecount2}). It is clear that

\be \label{crudebd} \sum_{k=2}^d (k-1) {\mathbf 1}_{\{[\eta_t] \in
\Net^{(k)} _{t-1}\}} \leq d-1 \, ,\ee hence (\ref{incdef}) implies
that $\E[ \xi_t \mid \F_{t-1}] \leq \eps$ and so $\E Y_t = O(\eps
t)$ and the process $\{\eps j - Y_j\}_{j\geq 0}$ is a
submartingale. Doob's maximal $L^2$ inequality (see \cite{D})
gives \be \label{lemdoob} \E[ \max _{j\leq t} (\eps j - Y_j)^2]
\leq 4\E [ (\eps t - Y_t)^2] \, .\ee

By (\ref{incdef}) and (\ref{uniform0}) we have
$$ \E[ \xi_j \mid \F_{j-1} ] \geq {(1+\eps) \over d-1} \cdot {d(d-1)
\NN_{j-1}^{(d)} - \AA_{j-1} \over dn -2j + 1} -1 \, .$$ By
(\ref{dneutral}) for all $j$ we have $\NN_{j-1}^{(d)} \geq n-2j$
and by (\ref{arec}) we have $\AA_{j-1} \leq d + (d-2)j$. We deduce
by the previous display that $ \E[\xi_j \mid \F_{j-1}] \geq
-D\eps$ for some fixed $D>0$ and all $j < C\eps n$. We learn that
for any $k < j < C\eps n$
$$ | \E [ \xi_j - D\eps | \F_{k} ] | = O(\eps) \, .$$
It follows that for any $k <j$
$$ \E[(\xi_j - D\eps)(\xi_k - D\eps)] = O(\eps^2) \, .$$
We deduce from the above that for $t < C\eps n$
$$\E [(\eps t - Y_t)^2] = 2 \sum _{j<k}^t \E[(\xi_j - \eps)(\xi_k -
\eps)]+ \sum _{j \leq t} \E [(\xi_j - \eps)^2] = O(\eps^2 t^2 + t)
\, .$$ By Jensen inequality and (\ref{lemdoob}) we get that
$$ \E [ \min _{j\leq t} (Y_j-\eps j)] = O(\eps t + \sqrt{t}) \,
,$$ and inequality (\ref{activecount2}) concludes the proof of
(\ref{aest}). As $Z_t = A_t - Y_t$, and $\E Y_t = O(\eps t)$ this
also concludes the proof of (\ref{zest}). \qed \\


\noindent {\bf Proof of Lemma \ref{moments}.} As $[w_t] \in
\Ne^{(d)}_{t-1}$ implies that $A_{t-1}=0$ we have by
(\ref{dneutral}) that
$$ \NN^{(d)}_t \geq \NN^{(d)}_{t-1} - {\bf 1}_{\{ [\eta_t] \in \Net^{(d)}
_{t-1}\}} - {\bf 1}_{ \{A_{t-1}=0\}} \, .$$ As $\NN^{(d)}_0 = n-1$
we learn that $\NN^{(d)}_t \geq n - t - \sum_{i=1}^{t-1} {\bf
1}_{\{A_i = 0\}}$ and so by the definition of $Z_t$ we have that
$\NN^{(d)}_t \geq n - t - Z_t$. Thus (\ref{zest}) of Lemma
\ref{activeest} gives that
$$ \E \NN^{(d)}_t \geq n - t - O(\eps t + \sqrt{t}) \, .$$
Also, by (\ref{dneutral}) we have that
$$ \NN^{(d)}_t \leq  \NN^{(d)}_{t-1} - {\bf 1}_{\{ [\eta_t] \in \Net^{(d)}
_{t-1}\}} \, .$$ Hence, (\ref{uniform0}) and $1-x \leq e^{-x}$
give that
$$ \E \Big [ \NN_t^{(d)} \mid {\cal F}_{t-1} \Big ] \leq \NN_{t-1}^{(d)}\Big ( 1  - {d \over dn-2t+1} \Big ) \leq
\NN_{t-1}^{(d)}e^{-{d \over dn-2t+1}} \, .$$  By iterating this we
get that
$$\E \NN_t^{(d)} \leq n e
^{-d \sum _{i=0}^t {1 \over dn-2i+1}} \leq ne^{-{t \over n}} \leq
n - t +{t^2 \over 2n} \, ,$$ where the last inequality is due to
$e^{-x} \leq 1 - x + x^2/2$ for all $x>0$. This concludes the
proof of (\ref{dbound}) as ${t^2 \over n} = O(\eps t)$ for $t <
C\eps n$.

Observe that (\ref{kneutral}) and (\ref{kneutral2}) implies that
$$ \E \Big [ \NN_{t}^{(d-1)} \mid {\cal F}_{t-1} \Big ] \leq \NN_{t-1}^{(d-1)}
+ (1-p) \, ,$$ which by iterating yields $\E N_t^{(d-1)} \leq
(1-p)t$. To complement this with a lower bound we use
(\ref{kneutral}) and (\ref{uniform0}) to get
$$ \E \Big [ \NN_{t}^{(d-1)} \mid {\cal F}_{t-1} \Big ]
\geq \NN_{t-1}^{(d-1)}+ (1-p){d\NN_{t-1}^{(d)} \over dn-2t+1}-
{(d-1)\NN_{t-1}^{(d-1)} \over dn-2t+1} - {\bf 1}_{\{A_{t-1}=0\}}
\, .$$ We now take expectation and bound the second term of the
right hand side using (\ref{dbound}) and the third term by
$\NN_{t-1}^{(d-1)} \leq t$ for all $t$. This yields
$$\E \NN_{t}^{(d-1)} \geq \E \NN_{t-1}^{(d-1)} + (1-p){d(n-t-O(\eps t + \sqrt{t})) \over dn-2t+1}
- O(t/n) - \prob (A_{t-1}=0) \, .$$ By iterating and using
(\ref{zest}) we get
$$ \E \Big [ N_{t}^{(d-1)} \Big ] \geq (1-p) \sum _{i=1}^t \Big (
1 - {(d-2)i+1+O(\eps i + \sqrt{i}) \over dn-2i+1} \Big ) -
O(t^2/n) - O(\eps t + \sqrt{t}) \, .$$ The sum can be bounded
below by $t- O(t^2/n)$ and as $t^2/n = O(\eps t)$ for $t \leq
C\eps n$ we conclude the proof of (\ref{d-bound}).

To prove the bound (\ref{kbound}) note that by (\ref{kneutral}) we
have
$$ \E [ \NN_t^{(k)} \mid \F_{t-1} ] \leq \NN_{t-1}^{(k)} + \prob \Big ( [\eta_t] \in
\Net^{(k+1)} _{t-1} \Big ) \, .$$ As $\NN_{t}^{(k+1)}\leq t$ for
$k<d-1$, using (\ref{uniform}) and iterating gives
(\ref{kbound}). \qed \\

\noindent {\bf Proof of Lemma \ref{deviation}.} By Lemma
\ref{moments} and the triangle inequality, the assertion of the
lemma is trivial for $k \in \{1, \ldots, d-2\}$ as ${t^2 \over n}
= O(\eps t)$ for $t < C\eps n$. We first prove the assertion for
$k=0$. By iterating (\ref{0neutral}) and (\ref{0neutral2}) we get
that

\begin{eqnarray*} \NN_t^{(0)} = 1 &+& \sum _{i=1}^t {\mathbf 1} _{\{(w_i, \eta_i) \in G_p\}} {\mathbf
1}_{\{\eta_i \hbox{ is neutral}\}} \nonumber \\
&+& \sum _{i=1}^t {\mathbf 1} _{\{(w_i, \eta_i) \not \in G_p\}}
{\mathbf 1}_{\{[\eta_i] \in \Net_{i-1}^{(1)} \}} +\sum _{i=1}^t
{\bf 1}_{\{w_{i+1} \hbox{ is neutral} \}}.
\end{eqnarray*}
Write
\begin{eqnarray*} &\,& X_1(t) = \sum _{i=1}^t {\mathbf 1} _{\{(w_i, \eta_i) \in
G_p\}} {\mathbf 1}_{\{\eta_i \hbox{ is neutral}\}} \, ,\\
&\,& X_2(t) = \sum _{i=1}^t {\mathbf 1} _{\{(w_i, \eta_i) \not \in
G_p\}} {\mathbf 1}_{\{[\eta_i] \in \Net_{i-1}^{(1)} \}} \, , \\
&\,& X_3(t) = \sum _{i=1}^t {\bf 1}_{\{w_{i+1} \hbox{ is neutral}
\}} \, .\end{eqnarray*} By definition $X_3(t) \leq Z_{t+1}$, hence
the triangle inequality implies

\be \label{devstep1} \E \Big | X_3(t) - \E X_3(t) \Big | \leq 2\E
Z_{t+1} = O(\eps t + \sqrt{t}) \, ,\ee where the last inequality
is due to (\ref{zest}). By (\ref{uniform}) and (\ref{kbound}) we
have for $i<C \eps n$
$$ \E \Big | {\mathbf 1}_{\{[\eta_i] \in {\bf N}_{i-1}^{(1)} \}} - \E {\mathbf
1}_{\{[\eta_i] \in {\bf N}_{i-1}^{(1)} \}} \Big | \leq {\eps i
\over n} \, ,$$ and hence the triangle inequality gives that

\be \label{devstep2} \E \Big | X_2(t) - \E X_2(t) \Big | \leq
O\Big ({\eps t^2 \over n} \Big ) = O(\eps t)\, . \ee By writing
${\mathbf 1}_{\{\eta_i \hbox{ is neutral}\}} = 1 - {\bf
1}_{\{\eta_i \in \AA_{i-1}\}}$ we get by the triangle inequality

\begin{eqnarray} \label{devstep3}
\E \Big | X_1(t) - \E X_1(t) \Big | &\leq&  \E \Big | \sum
_{i=1}^t {\mathbf 1} _{\{(w_i, \eta_i) \in G_p\}} - pt \Big | \\
&+& \E \Big | \sum_{i=1}^t {\mathbf 1} _{\{(w_i, \eta_i) \in
G_p\}}{\bf 1}_{\{\eta_i \in \AA_{i-1}\}} - p \sum_{i=1}^t \prob (
\eta_i \in \AA_{i-1}) \Big | \, . \nonumber
\end{eqnarray}
Since $\sum _{i=1}^t {\mathbf 1} _{\{(w_i, \eta_i) \in G_p\}}$ is
distributed as Bin$(t,p)$, the first expectation on the right hand
side of (\ref{devstep3}) is $O(\sqrt{t})$. By (\ref{uniform}) and
(\ref{aest}) of Lemma \ref{activeest} we get for each $i \leq t <
C \eps n$,
$$ \prob \Big ( \eta_i \in \AA_{i-1} \Big ) = {\E \AA_{i-1} \over
dn -2i +1} \leq O\Big ({\eps t + \sqrt{t} \over n}\Big ) \, .$$
Therefore, \be \E \Big | X_1(t) - \E X_1(t) \Big | \leq O(\eps t +
\sqrt{t}) \, . \nonumber \ee This together with (\ref{devstep1})
and (\ref{devstep2}) implies that
$$ \E | \NN_t^{(0)} - \E \NN_t^{(0)}| = O(\eps t + \sqrt{t}) \, .$$

We now prove the assertion of the lemma for $k=d-1$. After
choosing $w_{t+1}$ and before choosing $\eta_{t+1}$ we have $2t$
explored vertices and $\AA_t$ active vertices which belong only to
$d$-tuples with at most $d-1$ neutral vertices in them. Therefore,
$$ \AA_t + 2t = \sum _{k=0}^{d-1} (d-k) \NN_t ^{(k)} \, ,$$
and thus
$$ \NN_t^{(d-1)} = \AA_t + 2t - \sum _{k=0}^{d-2} (d-k) \NN_t^{(k)} \,
.$$ Hence the triangle inequality implies that,
$$ \E | \NN_t^{(d-1)} - \E \NN_t^{(d-1)}| \leq \E| \AA_t - \E \AA_t| + d \sum
_{k=0}^{d-2} \E |\NN_t^{(k)} - \E \NN_t^{(k)}| \, .$$ As we
verified the assertion of the lemma for $k\leq d-2$, by
(\ref{aest}) of Lemma \ref{activeest} we get the lemma for
$k=d-1$. The assertion for $k=d$ follows immediately by the
triangle inequality and the fact that
$$ \NN_t ^{(d)} = n - \sum_{k=0}^{d-1} \NN_t^{(k)} \, . $$ \qed \\

\noindent {\bf Proof of Corollary \ref{drift}.} We simply use
(\ref{uniform}) to plug into (\ref{incdef}) the bounds obtained in
Lemma \ref{moments}. We get

\begin{eqnarray*}
\E \xi_t &\leq& {1+\eps \over d-1} \Big [ {(d-1)d(n-t) \over
dn-2t+1} + {(d-1)(d-2)(1-p)t \over dn-2t+1} \Big ] - O\Big ( {\eps
t + \sqrt{t} \over n} \Big ) - 1 \, .
\end{eqnarray*}
Writing $1-p = {d-2-\eps \over d-1}$ and expanding the right hand
side gives that
$$ \E \xi _t - \eps + {d-2 \over d(d-1)} \cdot {t \over n}
= O \Big ( {\eps t + \sqrt{t} \over n} \Big ) = O\Big( \eps^2 +
{\sqrt{t} \over n} \Big ) \, ,$$ as $t \leq C\eps n$. This proves
part (i) of the corollary. Part (ii) follows immediately from
(\ref{uniform0}), Lemma \ref{activeest} and Lemma \ref{deviation}.
To prove part (iii), the bound on $\E[ \xi_t^2 \mid \F_{t-1} ]$,
we square (\ref{incdef}) and estimate it using Lemma \ref{moments}
and Lemma \ref{activeest}. For any $i\neq j$ we have ${\bf
1}_{\{[\eta_t] \in \Net^{(i)} _{t-1}\}} {\bf 1}_{\{[\eta_t] \in
\Net^{(j)} _{t-1}\}} = 0$, and also ${\bf 1}_{\{[\eta_t] \in
\Net^{(i)} _{t-1}\}} {\mathbf 1}_{\{ [\eta_t] \in \Act_{t-1}\}} =
0$. So by (\ref{incdef}) we have
$$ \E [ (\xi_t+1)^2 \mid \F_{t-1} ] = \prob ( (w_t, \eta_t) \in G_p )  \sum_{k=2}^d
(k-1)^2 \prob ( [\eta_t] \in \Net^{(i)} _{t-1}) - \prob (\eta_t
\in \Act_{t-1}) \, .$$

For any $k<d$, as $\NN_t^{(k)} \leq t$, by (\ref{uniform}) we have
$\prob \Big ([\eta_t] \in \Net^{(k)} _{t-1} \Big  ) = O(t/n)$ and
as $\AA_t \leq d+(d-2)t$ we have $\prob \Big (\eta _t \in
\Act_{t-1}\}\Big ) = O(t/n)$. As $\NN_{t-1}^{(d)} \geq n-2t$ by
(\ref{uniform0}) we have that $\prob \Big  ([\eta_t] \in
\Net^{(d)} _{t-1} \mid \F_{t-1} \Big ) = 1- O(t/n)$. All this
gives that
$$  \E [ (\xi_t+1)^2 \mid \F_{t-1} ] = p(d-1)^2(1- O(t/n)) -
O(t/n)= d - 1 + O(\eps) \, , $$ and as $\E [\xi_t \mid \F_{t-1} ]
= O(\eps)$ for $t< C\eps n$ we deduce that $\E [\xi_t^2 \mid \F_{t-1}]  = d-2 + O(\eps)$.  \qed \\

\noindent {\bf Proof of Lemma \ref{extradev}.} Note that for any
$t < \delta n$ we have $\NN_{t-1}^{(d)} \geq n(1-2\delta)$. Thus,
for such times $\NN_t^{(d)}$ can be stochastically bounded above
by $n - \sum_{j=1}^t I_j$ where $\{I_j\}$ are i.i.d. Bernoulli
random variables receiving $1$ with probability $1-2\delta$ and
$0$ with probability $2\delta$. By Large Deviation (see \cite{AS}
section A.14) we get (\ref{deviation2}).

By the same reasoning, for all times $t<\delta n$ the random
variable can be stochastically bounded below by $\sum_{i=1}^t J_i$
where $\{J_i\}$ are i.i.d Bernoulli random variables receiving $1$
with probability $p(1-2\delta)$ and $0$ with probability
$1-p(1-2\delta)$, which by Large Deviation yields (\ref{deviation3}). \qed \\

\noindent {\bf Proof of Lemma \ref{excessbd}.} We know that
$\NN_{t-1}^{(k)} \leq 2t$ for all $1 \leq k \leq d-1$ and that
$A_t \leq d+(d-2)t$ for all $t$. Thus by (\ref{uniform0}) for all
$t<T<n/4$ we have

$$ \prob ( \eta_t \in \cup _{k=1}^{d-1} \Net_{t-1}^{(k)} \cup \Act_{t-1} \mid \F_{t-1} ) \leq {(3d-4) t +d
\over dn-2t+1} \leq {4T \over n} \, .$$

Thus we can stochastically bound $|U_T|$ from above by a random
variable distributed as Bin$(T,q)$, where $q = { 4T \over n }$.
Thus, standard large deviations bounds, see Corollary $A.1.10$ of \cite{AS}, conclude the proof. \qed \\

\section{Inside the scaling window}

In this Section we prove Theorem \ref{criticalupper}. We follow
the strategy laid out in \cite{NP}. \\


\noindent {\bf Proof of Theorem \ref{criticalupper},
(\ref{part1})}. Let $\eps(n) = \lambda n^{-1/3}$ and $p={1 +
\eps(n) \over d-1}$. Let $\alpha$ be a random variable which
receives $d-2$ with probability $p$ and $-1$ with probability
$1-p$. Let $\{\alpha_i\}$ be i.i.d. random variables distributed
as $\alpha$ and let $\{W_t\}$ be the process defined by $W_t = d +
\sum_{i=1}^t \alpha_i$. By (\ref{incdef}), we can couple $\{Y_t\}$
and $\{W_t\}$ such that $Y_t \leq W_t$ for all $t$. Let $h =
n^{1/3}$ and define $\gamma = \gamma_h$ by
$$ \gamma = \min \{ t : W_t = 0 \textrm{ or } W_t \geq h \} \, .$$
For any $c>0$ we have
$$ \E e^{-c\alpha} = e^c\Big(1+p(e^{-c(d-1)}-1) \Big) \, ,$$
and by expanding both exponentials we get
$$ \E e^{-c\alpha} = (1+c+c^2/2+\ldots)\Big [ 1 + p(-c(d-1) +
c^2(d-1)^2 /2 - \ldots ) \Big ] \, .$$ It is straight forward to
check if we set $c=4\eps$, as long as $\eps>0$ is small enough, we
have $\E e^{-c \alpha} \geq 1$. Similarly, if $\eps < 0$ with
$|\eps|$ small enough we have that $\E e^{c\alpha} \leq 1$ for
$c=4\eps$. Thus, if $\lambda > 0$, part (i) of Lemma \ref{gambler}
and $1-e^{-x}\leq x$ for $x>0$ implies that

\be \label{gamb1} \prob ( W_\gamma > 0) \leq {4d\lambda \over
1-e^{-4\lambda}} n^{-1/3} \, .\ee A similar computation and an
application of part (ii) of Lemma \ref{gambler} shows that for
$\lambda < 0$ and $n$ large enough we have

\be \label{gamb2} \prob ( W_\gamma
> 0) \leq {-5d\lambda \over e^{-4\lambda}-1} n^{-1/3} \, .\ee
Also, when $\lambda = 0$ the process $\{W_t\}$ is a martingale and
we deduce by optional stopping that $\prob( W_\gamma > 0) \leq
dn^{-1/3}$. We now estimate $\E \gamma$ for all $\lambda$. Assume
first $\lambda
> 1/4$; as $\{ W_t - t \lambda n^{-1/3}\}$ is a martingale, the
optional stopping theorem gives
$$d = \prob(W_\gamma >0) \E [ W_\gamma \mid W_\gamma > 0 ] - \lambda n^{-1/3}\E\gamma \, .$$
We use $1-e^{-4\lambda} > 1/2$ for $\lambda > 1/4$ in
(\ref{gamb1}) and the fact that $\E [W_\gamma \mid W_\gamma >0]
\leq n^{1/3} + d$ to rearrange the last display. This gives that
$\E \gamma \leq 8dn^{1/3}$, for $\lambda > 1/4$. It is straight
forward to check that $\{W_t^2 - {1 \over 2}t\}$ is a
submartingale for any $\lambda > 0$, hence by optional stopping,
$$ 1 \leq \prob(W_\gamma >0) \E [W_\gamma^2 \mid W_\gamma > 0] - {1\over 2} \E \gamma \, .$$
We use ${4\lambda \over 1- e^{-4\lambda}} \leq 2$ for $\lambda \in
(0,1/4]$ in (\ref{gamb1}) and the obvious estimate $\E [W_\gamma^2
\mid W_\gamma > 0] \leq (n^{1/3}+d)^2$ to rearrange the last
display. This gives that $\E \gamma \leq 8dn^{1/3}$, for $\lambda
\in (0,1/4]$. An almost identical computation for the case
$\lambda \leq 0$ yields that for all $\lambda \in \mathbb{R}$ we
have
$$ \E \gamma \leq 8dn^{1/3} \, .$$

Define $\gamma^* = \gamma \wedge n^{2/3};$ by the last display,
inequalities (\ref{gamb1}) and (\ref{gamb2}) we deduce that there
exists $C=C(\lambda)$ such that

\be \label{midstep} \prob (W_\gamma^* > 0) \leq \prob (W_\gamma
\geq n^{1/3}) + \prob(\gamma \geq n^{2/3}) \leq Cn^{-1/3} \, . \ee

\noindent Taking an exponential in (\ref{incdef}) gives
$$ \E [ e^{c\xi_t} \mid \F_{t-1} ] = e^{-c} \E \Big [ e^{c{\mathbf 1} _{\{(w_t, \eta_t) \in G_p\}}
\Big ( \sum _{k=2}^d (k-1) {\mathbf 1}_{\{[\eta_t] \in \Net^{(k)}
_{t-1}\}} - {\mathbf 1}_{\{ \eta_t \in \Act_{t-1}\}}  \Big )} \mid
\F_{t-1}\Big ] \, .$$ The conditional expectation on the right
hand side of the last display is $e^{c(k-1)}$ with probability
${pk\NN_{t-1}^{(k)} \over dn-2t+1}$ for any $2 \leq k \leq d$ by
(\ref{uniform0}) and at most $1$ with probability $1-p+{p\AA_{t-1}
\over dn-2t+1}$ by \ref{uniform2}. Thus,
$$ \E [ e^{c\xi_t} \mid \F_{t-1} ] \leq e^{-c} \Big [1+p\Big (-1+
\sum_{k=1}^{d} e^{c(k-1)}{k \NN_{t-1}^{(k)} \over dn -2t + 1} +
{\AA_{t-1} \over dn-2t+1} \Big ) \Big ] \, .$$ Using $e^x \leq
1+x+x^2$ for $x\in[0,1]$ we have that for $c<{1 \over d-1}$
\begin{eqnarray*} \E [ e^{c\xi_t} \mid \F_{t-1} ] \leq e^{-c} \Big [1+p \Big
(&-&1 + \sum _{k=1}^d \Big (1 + c(k-1) + {c^2(k-1)^2} \Big )
{k\NN_{t-1}^{(k)} \over dn - 2t+1} \\ &+& {\AA_{t-1} \over
dn-2t+1} \Big ) \Big ] \, .\end{eqnarray*} We expand the right
hand side of the last display using the fact that ${\AA_{t-1} +
\sum _{k=1}^d k \NN_{t-1}^{(k)} \over dn -2t +1} = 1$ and that
$\sum _{k=1}^{d-1} k \NN_{t-1}^{(k)} \leq (d-1)(n- \NN_{t-1}^{(d)}
- \NN_{t-1}^{(0)})$. This gives,
\begin{eqnarray} \label{insidestep} \E [ e^{c\xi_t} \mid \F_{t-1} ] \leq e^{-c} \Big [1+p
&\Big(& {c(d-1)dN_{t-1}^{(d)} + c(d-2)(d-1)(n- N_{t-1}^{(d)} -
N_{t-1}^{(0)}) \over dn-2t+1} \nonumber \\ &\,& + c^2 (d-1)^2 \,\,
\Big ) \Big ] \, .\end{eqnarray} For some small $\delta>0$ denote
by $\A$ the event
$$\A = \{ N_t^{(d)} \leq n - (1-\delta)t \, , \quad N_t^{(0)} >
(1-\delta)pt \, , \quad \forall n^{1/3}<t< \delta n /3 \}\, .$$ We
now condition on $\A$ and put $p={1+\eps \over d-1}$ in
(\ref{insidestep}). A straightforward computation yields that for
$c<{1 \over d-1}$ and $n^{1/3}<t< \delta n /3$ we have
$$ \E [ e^{c\xi_t} \mid \F_{t-1}, \A ] \leq e^{-c} \Big [1+(1+\eps) \Big
( c - {(d-2 + O(\delta))ct \over d(d-1)n} +c^2(d-1) \Big ) \Big ]
\, .$$ As $d\geq 3$ we can choose $\delta$ small enough such that
${(d-2 + O(\delta)) \over d} > {1 \over 4}$; we also use $1+x \leq
e^x$ for all $x>0$ in the last display. This gives that for such
$t$'s,

\be \label{computation0} \E [ e^{c\xi_t} \mid \F_{t-1}, \A ] \leq
e^{c\eps -{ct \over 4(d-1)n} +2(d-1)c^2} \, .\ee By estimating
$e^{c\xi_j} \leq e^{c(d-2)}$ for all $j\leq n^{1/3}$, as $\gamma*
\leq n^{2/3}$ we get from the last display that for any $t< \delta
n/3 - n^{2/3}$

\be \label{computation}  \E \Big [ e^{ c \sum _{j=1}^{t}
\xi_{\gamma^*+j}} \mid \gamma^* , \, \A \Big ] \leq e^{c\eps t -
{c t^2 \over 8(d-1)n} +2(d-1)c^2 t + c(d-2)n^{1/3} } \, . \ee

Define the process $\{R_t\}$ by
$$ R_t = Y_{\gamma^*+t} - Y_{\gamma^*} = \sum _{j=1}^t \xi_{\gamma^*+j} \, .$$
As the estimate (\ref{computation}) is uniform in $W_{\gamma^*}$
and $\gamma^*$ we get that
$$ \E [ e^{cR_t} \mid W_{\gamma^*}, \A ] \leq e^{c\eps t - {c t^2 \over 8(d-1)n} +2(d-1)c^2t +
c(d-2)n^{1/3} } \, .$$ Write $\prob_W$ for the conditional
probability measure given $W_{\gamma^*}$ and $\A$. Then by
previous equation, for any $c<{1 \over d-1}$ and $t<\delta n
/3-n^{2/3}$ we have
\begin{eqnarray*}
\prob_W \Big ( R_t \geq -W_{\gamma^*} \Big ) &\leq& \prob_W \Big (
e^{cR_t} \geq e^{-cW_{\gamma^*}} \Big ) \\ &\leq&  e^{c\eps t - {c
t^2 \over 8(d-1)n} +2(d-1)c^2t + c(d-2)n^{1/3} }e^{cW_{\gamma^*}}
\, .
\end{eqnarray*}
By (\ref{deviation2}) and (\ref{deviation3}) of Lemma
\ref{extradev} it follows that $\prob(\A ^c) \leq n e^{-a
n^{1/3}}$ for some fixed $a>0$. As $Y_{\gamma^*} \leq
W_{\gamma^*}$ it follows by the definition of $R_t$ and by
conditioning on $\A$ that
\begin{eqnarray*} \prob \Big ( Y_{\gamma^*+t} > 0 \mid W_{\gamma^*}>0 \Big )
\leq \E \Big [ \prob_W  (R_t \geq -W_{\gamma^*}) \mid W_{\gamma^*}
> 0 \Big ] + \prob( \A ^c) \\ \leq e^{c\eps t - {c t^2 \over
8(d-1)n} +2(d-1)c^2t + c(d-2)n^{1/3} } \E \Big [e^{cW_{\gamma^*}}
\mid W_{\gamma^*} > 0 \Big ] + n e^{-a n^{1/3}} \, .
\end{eqnarray*}
Since $W_{\gamma^*} \leq n^{1/3}+d$ we can bound the conditional
expectation on the right hand side. This yields,
$$ \prob \Big ( Y_{\gamma^*+t} > 0 \mid W_{\gamma^*}>0 \Big ) \leq
e^{c\eps t - {c t^2 \over 8(d-1)n} +2(d-1)c^2t + c(d-1)n^{1/3} +
cd } + n e^{-a n^{1/3}} \, .$$ Now recall that $\eps = \lambda
n^{-1/3}$ and take $c = {{t^2 \over 8(d-1)n} - \eps t -(d-1)
n^{1/3} \over 4t(d-1)}$ and $t=Bn^{2/3}$ for some $B>0$ large
enough so that $c>0$. Note that $c$ is of order $n^{-1/3}$ so
clearly $c < {1 \over d-1}$. Putting all this together gives that

\begin{eqnarray*} \prob \Big ( Y_{\gamma^*+Bn^{2/3}} > 0 \mid W_{\gamma^*}>0
\Big ) &\leq& e^{ - {\Big ({B^2 \over 8(d-1)} - \lambda B
-(d-1)\Big )^2 \over 8B(d-1)} + O(n^{-1/3})} + ne^{-\alpha
n^{2/3}} \\ &\leq& e^{-r B^3} \, ,\end{eqnarray*} for some
$r=r(\lambda)>0$ and $n$ large enough. Recall that $t_1$ is the
first time the process $Y_t$ hits $0$ and that Lemma
\ref{bottomline} implies that $|\C(v)| \leq t_1$. Thus, by our
coupling, if $|\C(v)|>An^{2/3}$ then $W_{\gamma^*}>0$ and
$Y_{\gamma^* + (A-1)n^{2/3}} > 0$. Thus by the previous inequality
and (\ref{midstep}), for $A>1$ we have
$$ \prob(|\C(v)|\geq An^{2/3}) \leq Cn^{-1/3}e^{-r(A-1)^3} \, .$$
Denote by $N_{T}$ the number of vertices contained in components
larger than $T$. Observe that $|\C_1| \geq T$ implies $N_T \geq
T$. So taking $T=An^{2/3}$ gives
\begin{eqnarray*} \prob \Big ( |{\cal C}_1| \geq T \Big ) &\leq& \prob
\Big ( N_{T} \geq T \Big ) \leq \frac{\E N_{T}}{T}
 \leq \frac{n \prob(|\C(v)|\geq T)}{T} \leq {C\over A}e^{-r(A-1)^3} \, ,
\end{eqnarray*} concluding the proof. \qed \\


\noindent {\bf Proof of Theorem \ref{criticalupper},
(\ref{part2})}. Let $\delta \in (0,1)$ be small and let
$\gamma=\gamma(\delta,\lambda)>0$ be determined later. Put
$h=\gamma n^{1/3}$, $T_1=n^{2/3}$ and $T_2=\delta n^{2/3}$. As in
\cite{NP} we ensure that with high probability the process
$\{Y_t\}$ gets to height $h$ before time $T_1$, and then stays
positive for at least $T_2$ steps. This ensures by Lemma
\ref{bottomline} that $|{\cal C}_1| > { \delta n^{2/3} \over d-1}$
with high probability. Indeed, let us define the stopping time
$$ \tau _h = \min \{ t \leq T_1 : Y_t \geq h\} \,$$
if this set is nonempty, and $\tau _h = T_1$ otherwise. Observe
that $Y_t^2 - Y_{t-1}^2 = \xi_t^2 + 2\xi_t Y_{t-1}$. By Corollary
\ref{drift}, we have $\E [\xi_t^2 \mid \F_{t-1}] = d-2 +
O(n^{-1/3})$ and $\E [ \xi_t \mid \F_{t-1} ] \geq
\Omega(n^{-1/3})$ for $t \leq T_1$. Thus, if $Y_{t-1} \leq h$ and
$\gamma$ is small enough, we have that for all $t\leq T_1$,
$$ \E \Big [ Y_t^2 - Y_{t-1}^2 \, \Big | \, Y_{t-1} \Big ] \geq {1
\over 2} \, .$$ Hence $Y_{t \wedge \tau_h} ^2 - (t \wedge \tau_h)
/2$ is a submartingale. As $Y_{\tau_h} \leq h+d$ we have that $\E
Y_{\tau_h}^2 \leq (h+d)^2 \leq 2h^2$, so by optional stopping we
get
$$ 2h^2 \geq \E Y_{\tau_h}^2 \geq {1 \over 2} \E
\tau _h \geq {T_1 \over 2} \prob \Big ( \tau_h = T_1 \Big ) \, ,$$
hence

\be \label{top} \prob \Big ( \tau_h = T_1 \Big ) \leq {4h^2 \over
T_1 } \, . \ee

Write $\prob_h$ for conditional probability given the event $\{
\tau_h < T_1 \}$ and $\E_h$ for conditional expectation given that
event and define
$$ \tau _0 = \min \{t \leq T_2 : Y_{\tau_h + t} = 0\} \, .$$
We wish to bound from above the probability that $\tau_0 \leq T_2$
given that $\tau_h < T_1$. As before there exists a constant
$C=C(\lambda)$ such that $\E[ \xi_t \mid \F_{t-1} ] \geq
-Cn^{-1/3}$ for all $t\leq T_1 + T_2$. Thus the process
$$ S_t = Y_{\tau_h +t} - Y_{\tau_h} + tCn^{-1/3} \, $$
is a submartingale and hence so is $S_t^2$. We conclude that as
long as $h > T_2 C n^{-1/3}$
\begin{eqnarray} \label{doobstep}
\prob_h \Big ( \min _{t \leq T_2} Y_{\tau_h +t} \leq 0 \Big )
&\leq& \prob_h \Big ( \min _{t \leq T_2} S_t \leq -h +
T_2Cn^{-1/3} \Big ) \nonumber \\ &\leq& \prob_h \Big ( \max _{t
\leq T_2} S_t^2 > (h-T_2Cn^{-1/3})^2 \Big ) \nonumber \\ &\leq&
{4\E_h S_{T_2}^2 \over (h-T_2Cn^{-1/3})^2 } \, ,
\end{eqnarray}
where the last inequality is Doob's Maximal inequality (see
\cite{D}). As usual, for any $k<j<T_1+T_2$ we can bound
$$ \E [ \xi_j - Cn^{-1/3} | \F_{k-1} ] = O(n^{-1/3}) \, ,$$
and so
$$ \E[(\xi_j - Cn^{-1/3})(\xi_k - Cn^{-1/3})] = O(n^{-2/3}) \, .$$
This together with the fact that $\xi_{\tau_h +j} - Cn^{-1/3}$ is
bounded by $d-2$ shows that
\begin{eqnarray*} \E_h [ S_{T_2}^2 \mid \tau_h ] &=& \sum _{j \neq
k}^{T_2} \E_h [ (\xi_{\tau_h+j} - Cn^{-1/3})(\xi_{\tau_h+k} -
Cn^{-1/3}) \mid \tau_h ] \\ &+& \sum _{j=1}^{T_2} \E_h [
(\xi_{\tau_h+j} - Cn^{-1/3})^2 \mid \tau_h ] = O(n^{-2/3}T_2^2 +
T_2) = O(\delta n^{2/3}).
\end{eqnarray*}
Hence (\ref{doobstep}) implies that
$$\prob _h (\tau_0 \leq T_2) \leq { O(\delta n^{2/3}) \over (h-\delta C n^{1/3})^2 } \, ,$$
as long as $\gamma > \delta C$, so the denominator is positive.
Combining this with (\ref{top}) gives
\begin{eqnarray*} \prob (\tau _0 \leq T_2) \leq \prob (\tau_h = T_1) +
\prob _h (\tau_ 0 \leq T_2) &\leq& {4h^2 \over T_1 } + { O(\delta n^{2/3}) \over (h-\delta C n^{1/3})^2 } \\
 &=& 4\gamma^2 + {O(\delta)
\over (\gamma - \delta C)^2} \, ,
\end{eqnarray*}
and by choosing $\gamma = \delta C + \delta^{1/4}$ we deduce that
$$ \prob (\tau _0 \leq T_2) \leq D\delta^{1/2} \, ,$$
for some constant $D=D(\lambda)>0$. By Lemma \ref{bottomline},
$|\C_1| \leq {T_2 \over d-1}$ implies
$\tau_0 \leq T_2$, which concludes the proof. \qed \\

\section{Below the scaling window}

We use Lemma \ref{tail} on another specific case. Fix some small
$\eps>0$ and set $p = {1 - \eps \over d-1}$. Let $\beta$ be a
random variable receiving $d-2$ with probability $p$ and $-1$ with
probability $1-p$. Let $\{W_t\}$ and $\tau$ be defined as in Lemma
\ref{tail} with $W_0=d$.

\begin{lemma} \label{hitting} There exists constant $c_1, c_2> 0$ such that for all $T>\eps^{-2}$ we have
$$ \prob (\tau \geq T ) \geq c_1\Big ( \eps^{-2} T^{-3/2} e^{- {(\eps^2 + c_2 \eps^3) T \over
2(d-2)}}\Big ) \, .$$ Furthermore,
$$ \E \tau^2 = O ( \eps^{-3} ) \, .$$
\end{lemma}

\noindent {\bf Proof.} We estimate $\theta_0$ defined in Lemma
\ref{tail}. By (\ref{tailassump1}) we have
$$ \varphi(\theta_0) ^{-1} \Big [ {e^{\theta_0 (d-2)}(d-2)(1-\eps) \over
d-1} - {e^{-\theta_0} (d-2+\eps) \over d-1} \Big ] = 0 \, .$$ By
estimating $e^{x} = 1+x+O(x^2)$ we get
$$ \theta_0 = {\eps \over d-2} + O(\eps^2) \, .$$ We have
$$ \varphi(\theta) = {1 - \eps \over d-1}  e^{\theta (d-2)} +
{d-2+\eps \over d-1} e^{-\theta} \, .$$ Plugging in the value of
$\theta_0$ and writing $e^x = 1 + x + {x^2 \over 2} + O(x^3)$
gives that
$$ \varphi(\theta_0) = 1 - {\eps ^2 \over 2(d-2)} + O(\eps^3) \,
.$$ Thus by Lemma \ref{tail} we have

\begin{eqnarray*}
\prob( \tau \geq T) = \sum _{\ell \geq T} \Theta \Big (
\ell^{-3/2} \varphi(\theta_0) ^\ell \Big ) \, .
\end{eqnarray*}
Using our estimate on $\varphi(\theta_0)$ and the assumption that
$T>\eps^{-2}$ an immediate computation yields the first assertion
of the lemma. The second assertion follows from the following
computation. By Lemma \ref{tail} we have

$$ \E \tau ^2 = \sum _{\ell \geq 1} \ell^2 \prob (\tau = \ell) =
\leq C \sum _{\ell \geq 1} \sqrt{\ell} \varphi(\theta_0)^\ell \,
.$$ Thus, by direct computation (or by \cite{F}, section XIII.5,
Theorem 5)
$$ \E \tau ^2 \leq O \Big ( {1 \over 1-\varphi(\theta_0)} \Big)^{3/2} =
O(\eps^{-3}) \, .$$ \qed \\




\noindent {\bf Proof of Theorem \ref{lower}.} Note that
Proposition \ref{dregular} proves the upper bound on $|\C _\ell|$
implied in Theorem \ref{lower}, so we only need to prove the lower
bound. Write
$$ T = {2(1-\eta)(d-2)} \eps^{-2}\log(n\eps^3) \, .$$
For each integer $j\geq 0$ let $\{ W^{(j)}_t\}$ be independent
processes defined by $W^{(j)}_0 = Y_{t_j}$ and $W^{(j)}_t -
W^{(j)}_{t-1}$ receives $d-2$ with probability ${1 - (1+{\eta
\over 4})\eps \over d-1}$ and $-1$ otherwise. Note that for each
$j$, the process $\{W^{(j)}_t\}$ is just the process defined in
Lemma \ref{hitting} with $p={1 - (1+{\eta \over 4})\eps \over
d-1}$. By (\ref{incdef}) and (\ref{uniform0}), the variable
$\xi_t$ can always be stochastically bounded below by a variable
taking the value $d-2$ with probability ${1-\eps \over
d-1}\cdot{d\NN_{t-1}^{(d)} \over dn}$ and $-1$ otherwise. Since
$\NN_t^{(d)} \geq n-2t$ for all $t$, as long as $t < {\eta \over
8}\eps n$ we can stochastically bound $\xi_t$ below by $W^{(j)}_t
- W^{(j)}_{t-1}$. Thus, as long as $t_{j+1} < {\eta \over 8}\eps
n$, we can couple $\{Y_t\}$ and $W^{(j)}_t$ such that $Y_{t_j + t}
\geq W^{(j)}_t$ for all $t \in [0, t_{j+1}-t_j]$. Define the
stopping times $\{\tau_j\}$ by
$$ \tau_j = \min \{ t : W^{(j)}_t = W^{(j)}_0 - 1 \} \, .$$
By our coupling, it is clear that if $\tau_j > T$ then $t_{j+1} -
t_j > T$. Take
$$ N = \Big \lfloor \eps^{-1} (n\eps^3)^{(1-{\eta
\over 8})} \Big \rfloor \, .$$ We will prove that with high
probability $t_N < {\eta \over 8} \eps n$ and that there exists
$k_1 < k_2 < \ldots < k_\ell < N$ such that $\tau _{k_i} > T$.
Since $\E [ \xi_t \mid \F_{t-1} ] \leq -\eps$, we have by optional
stopping that $\E[t_{j+1}-t_j] \leq d\eps^{-1}$, and hence $\E t_N
\leq d \eps^{-2} (n\eps^3)^{(1-{\eta \over 8})}$ which implies
that

\be \label{part1ofproof} \prob \Big ( t_N > {\eta \over 8} \eps n
\Big ) \leq {8(n\eps^3)^{-{\eta \over 8}}\over \eta}  \to 0 \,
.\ee Also, by Lemma \ref{hitting} we have for some $c>0$
$$ \prob ( \tau _j > T ) \geq c \eps (n\eps^3)^{-(1+{\eta \over 4})^2(1-\eta)(1-c_2 \eps)} \log(n\eps^3)^{-3/2} \geq \eps (n\eps^3)^{-(1-{\eta \over 4})}\, , $$
as long as $\eta < 4$ and $\eps$ is small enough. Let $X$ be the
number of $j \leq N$ such that $\tau_j > T$. Then we have
$$ \E X \geq N \eps (n \eps^3)^{-(1-{\eta \over 4})} \geq C
(n\eps^3)^{{\eta \over 8}} \to \infty \, ,$$ hence by Large
Deviations (see \cite{AS}, section A.14) for any fixed integer
$\ell>0$ we have for some $c>0$

\be \label{part2ofproof} \prob \Big ( X < \ell \Big ) \leq e^{-c
(n\eps^3)^{{\eta \over 8}}} \to 0 \, .\ee Our coupling and Lemma
\ref{bottomline} imply that
$$ \Big \{ |\C_\ell| < {T \over d-1} \Big \} \subset \Big \{ X < \ell \Big \}
\cup \Big \{ t_N > {\eta \over 8} \eps n \Big \} \, ,$$ and hence
by (\ref{part1ofproof}) and (\ref{part2ofproof}) we have
$$ \prob \Big ( |\C _\ell| < {T \over d-1} \Big ) \to 0 \, . $$
\qed \\

\section{Above the scaling window}



We split the proof of Theorem \ref{upper} into two steps. In the
first step we show there is a unique component of order ${2d \over
d-2} \eps n$ which has about $2d \eps n$ closed edges separating
it from its boundary. In the second step we condition on this
event and restart the exploration process on the graph remaining
after removing this partial matching to get the estimates on the
$\ell$-th largest component for $\ell \geq 2$. The first step
follows the strategy laid out in \cite{NP2}.

We require some definitions. Consider $p$-bond percolation on the
configuration model, i.e., we draw a perfect matching on the
vertex set $\{ \, v_{i,k} \, : \, 1 \leq i \leq n \, , 1 \leq k
\leq d \}$ and then retain each edge with probability $p$ and
delete it with probability $1-p$ independently of all other edges.
Denote the resulting graph by $M(n,d,p)$ and recall that
$G^*(n,d,p)$ is the graph obtained from $M(n,d,p)$ by contracting
every tuple to a vertex. A set of $d$-tuples $S$ in $M(n,d,p)$ is
called a {\em component} if the vertex set corresponding to $S$ in
$G^*(n,d,p)$ is a connected component. We say that a $d$-tuple $v$
is {\em $k$-damaged}, with $0 \leq k \leq d$,  by a component $S$
if $v \not \in S$ and there are precisely $k$ closed edges (i.e.,
edges not retained in percolation) between a vertex in $v$ and a
vertex in a tuple belonging to $S$. Let $M_k(S)$ be the set of all
$k$-damaged tuples of a component $S$. Let $p={1 +\eps(n) \over
d-1}$. We say a component $S$ is {\em $\delta$-giant} for some
$\delta>0$ if the following properties hold:
\begin{enumerate}
\item[(i)] $(1-\delta){2d \over d-2} \eps n \leq |S| \leq
(1+\delta){2d \over d-2} \eps n \, ,$

\item[(ii)] $(1-\delta)2d \eps n \leq |M_1(S)| \leq (1+\delta)2d\eps n
\, .$ \\
\end{enumerate}

For $\delta>0$ let $\G(\delta)$ denote the event that there exists
a unique {\em $\delta$-giant} component in $M(n,d,p)$. The
following Theorems imply Theorem \ref{upper}.

\begin{thm}\label{upperthm1} Let $\eps(n)>0$ be a sequence such that $\eps(n) \to 0$ and
$\eps(n) n^{1/3} \to \infty$. Let $p = {1+\eps(n) \over d-1}$ and
consider $M(n,d,p)$. Then for any $\delta>0$ we have \be
\label{abovestep1} \prob ( \G(\delta) ) \to 1 \, , \qquad
\hbox{\rm as } n \to \infty .\ee
\end{thm}

\begin{thm}\label{upperthm2} Condition on $\G(\delta)$ and denote by
$S_1$ the $\delta$-giant component. Let $\{S_\ell\}_{\ell \geq 2}$
denote the components of $M(n,d,p)$ after removing $S_1$, ordered
by size. Then under the conditions of the previous theorem, for
any $\eta>0$ there is $\delta>0$ small enough such that

\be \label{abovestep2} \prob \Big ( |S_2| \geq (1+\eta) {2(d-2)
\over d-1} \eps^{-2}(n)\log (n \eps^3(n)) \, \Big | \, \G(\delta)
\Big ) \to 0 \, ,\ee and for any fixed integer $\ell \geq 2$ we
have

\be \label{abovestep3} \prob \Big ( |S_\ell| \leq (1-\eta) {2(d-2)
\over d-1} \eps^{-2}(n)\log (n \eps^3(n)) \, \Big | \, \G(\delta)
\Big ) \to 0 \, .\ee \\
\end{thm}

\noindent {\bf Proof of Theorem \ref{upper}}. Fix some $\eta>0$
and take $\delta>0$ small enough guaranteed by Theorem
\ref{upperthm2}. Theorem \ref{upperthm1} guarantees that the event
$\G(\delta)$ holds with high probability. Hence with that
probability and hence there exists a component of size between
$(1-\delta){2d \over d-2} \eps n$ and $(1+\delta){2d \over d-2}
\eps n$. We condition on $\G(\delta)$ and remove this component;
Theorem \ref{upperthm2} then implies that with high probability
the graph remaining has no components of size bigger than
$(1+\eta) {2(d-2) \over d-1} \eps^{-2}(n)\log (n \eps^3(n))$ and
the $\ell$-th component is bigger than $(1-\eta) {2(d-2) \over
d-1} \eps^{-2}(n)\log (n \eps^3(n))$. As these probabilities tend
to $0$ in the space $G^*(n,d,p)$, as the event $\Simp$ has
positive probability, we conclude the same for the space
$G(n,d,p)$. \qed \\

\noindent {\bf Proof of Theorem \ref{upperthm1}.} Write
$T={(1+\delta)2d(d-1)\over d-2} \eps n$ and $\xi^*_j = \E[\xi_t
\mid \F_{t-1}]$. The process
$$M_t = Y_t - \sum _{j=1}^t \xi^*_j \, ,$$
is a martingale. By Doob's maximal $L^2$ inequality (see \cite{D})
we have
$$ \E (\max _{t \leq T} M_t)^2 \leq 4\E M_T^2 \, .$$
As $M_t$ has orthogonal bounded increments we conclude $\E M_T^2 =
O(T)$. By Jensen inequality

\be \label{upperstep1} \E \Big [ \max _{t\leq T} \Big ( Y_t - \sum
_{j=1}^t \xi^*_j \Big ) \Big ] \leq O(\sqrt{T}) = O ( \sqrt{\eps
n} ) \, . \ee By (\ref{incdef}) and (\ref{uniform0}) for any
$j\leq T$ we have
$$ \xi^*_j - \E\xi_j = {1+\eps \over d-1} \sum_{k=2}^d (i-1)\Big [
{k \NN_{j-1}^{(k)} - i \E \NN_{j-1}^{(k)} \over dn-2j+1} \Big ] -
{\AA_{j-1} - \E \AA_{j-1} \over dn-2j+1} \, .$$ Applying the
triangle inequality to the last display, together with
(\ref{aest}) of Lemma \ref{activeest} and Lemma \ref{deviation}
gives that $\E | \xi^*_j - \E \xi_j | \leq O \Big ({\eps j +
\sqrt{j} \over n} \Big )$. So for any $t \leq T$ we have
$$ \E \Big [ \sum _{j =1}^t  | \xi_j^* - \E \xi_j
 | \Big ] = O(\eps^3 n) \, .$$
By the triangle inequality we get \be \label{upperstep2} \E \Big [
\max _{t \leq T} \Big | \sum _{j=1}^t ( \xi_j^* - \E \xi_j) \Big |
\Big ]\leq O(\eps^3 n)  \, . \ee Using the triangle inequality
together with (\ref{upperstep1}), (\ref{upperstep2}) and Markov's
inequality gives \be \label{keystep} \prob \Big ( \max _{t\leq T}
\Big | Y_t - \sum _{j=1}^t \E \xi_j \Big | \geq \delta \eps^2 n
\Big ) \leq \delta^{-1} ( O(\eps) + O ( (\eps^3 n)^{-1/2} ))
\longrightarrow 0 \, . \ee By Corollary \ref{drift} we have that
for any $b>0$ \be \label{conc2} \sum _{t=1}^{b\eps n} \E \xi_t =
\Big ( b-{(d-2)b^2 \over 2d(d-1)} \Big ) \eps^2n + O(\eps^3 n) \,
.\ee

Write
$$t' = {2d(d-1)\over (d-2)} \eps n \, .$$
Inequalities (\ref{keystep}) and (\ref{conc2}) imply that for
small $\delta > 0$ with probability tending to $1$, we have that
$Y_t$ is positive at times $[\delta t'/2, t'(1-\delta/2)]$. This
together with Lemma \ref{bottomline} implies that with high
probability we have explored a component containing at least
$(1-\delta) {2 d \over d-2} \eps n$ tuples.
Furthermore, by (\ref{keystep}) and (\ref{conc2}) we infer that

$$\prob \Big ( Y_{t'(1+\delta)} \leq - {2d(d-1)\over d-2} \delta (1+\delta) \eps^2 n + O(\eps^3 n) \Big ) \to 1 \, ,$$
and
$$ \prob \Big ( \forall \, t \leq \delta t' /2 \quad  Y_t > O(-\eps^3 n) \Big ) \to 1 \, .$$
Thus, with high probability, by time $t'(1+\delta)$ we have
completely explored a component of size at least $(1-\delta) {2 d
\over d-2} \eps n$. On the other hand, Lemma \ref{excessbd} and
Lemma \ref{bottomline} show that with high probability the size of
this component is at most $(1+\delta){2 d \over d-2} \eps n$.
Denote this component by $S$. By (\ref{d-bound}) of Lemma
\ref{moments} and Lemma \ref{deviation} we have that with high
probability $| \NN_{t'}^{(d-1)} - 2d\eps n| \leq \delta \eps n$.
This implies that $|M_1(S) - 2d\eps n| \leq \delta \eps
n$ with high probability and concludes our proof. \qed \\

To prove Theorem \ref{upperthm2} we need the following lemma,
which is just another application of Lemma \ref{tail} to a
specific case. Fix some small $\eps>0$ and let $\beta$ be a random
variable taking the value $d-2$ with probability ${1 - (2d-3)\eps
\over d-1}$, the value $d-3$ with probability $2\eps$ and the
value $-1$ with probability ${d-2 - \eps \over d-1}$. Let
$\{W_t\}$ and $\tau$ be defined as in Lemma \ref{tail} with
$W_0=d$.

\begin{lemma} \label{hitting2} There exists constant $C_1, C_2, c_1, c_2> 0$ such that for all $T>\eps^{-2}$ we have
$$ \prob (\tau \geq T ) \leq C_1\Big ( \eps^{-2} T^{-3/2} e^{- {(\eps^2 - c_1 \eps^3) T \over
2(d-2)}}\Big ) \, ,$$ and
$$ \prob (\tau \geq T ) \geq c_1\Big ( \eps^{-2} T^{-3/2} e^{- {(\eps^2 + c_2 \eps^3) T \over
2(d-2)}}\Big ) \, .$$ Furthermore,
$$ \E \tau^2 = O ( \eps^{-3} ) \, .$$
\end{lemma}

\noindent {\bf Proof.} We estimate $\theta_0$ of Lemma \ref{tail}.
By (\ref{tailassump1}) we have
$$ \varphi(\theta_0) ^{-1} \Big [ {e^{\theta (d-2)}(d-2)(1-(2d-3)\eps) \over
d-1} + e^{\theta_0 (d-3)}(d-3)2\eps - {e^{-\theta_0} (d-2-\eps)
\over d-1} \Big ] = 0 \, .$$ By estimating $e^{x} = 1+x+O(x^2)$ we
get
$$ \theta_0 = {\eps \over d-2} + O(\eps^2) \, .$$ We have
$$ \varphi(\theta) = {1 - (2d-3) \eps \over d-1}  e^{\theta (d-2)}
+ 2\eps e^{\theta_0 (d-3)} + {d-2-\eps \over d-1} e^{-\theta} \,
.$$ By estimating $e^x=1+x+x^2/2 + O(x^3)$ and plugging in the
value of $\theta_0$, we obtain that
$$ \varphi(\theta_0) = 1 - {\eps ^2 \over 2(d-2)} + O(\eps^3) \,
.$$ The rest of the proof is identical to the proof of Lemma
\ref{hitting}. \qed \\

\noindent {\bf Proof of Theorem \ref{upperthm2}.} Let $S$ be the
component specified in the event $\G(\delta)$. Condition on
$\G(\delta)$ and consider the graph remaining after removing $S$.
Denote by $\prob_S$ denote the distribution of this remaining
graph conditioned on $S$ and on the edges in the matching adjacent
to vertices in the tuples of $S$. Denote by $\prob_M$ the
distribution of $p$-bond percolation on a uniform matching on a
set of $\sum_{k=1}^d M_k(S)$ tuples of which precisely $M_k(S)$
tuples are of size $d-k$. Observe that $\prob_S$ is just $\prob_M$
conditioned on the event that the resulting graph has no
$\delta$-giant component. Theorem \ref{upperthm1} guarantees that
with high probability there is a unique $\delta$-giant component.
We learn that for any set of graphs $\B$ which do not contain an
$\delta$-giant component we have $\prob_S(\B) =
(1+o(1))\prob_M(\B)$. Thus it suffices to prove the required tail
bounds on the components in $\prob_M$. We do this in a similar
manner to the proof Theorem \ref{lower}.

Given $S$, the exploration process on the remaining graph,
starting from a tuple $v$ has the same dynamics described in
Section \ref{explore}. As $S$ is a $\delta$-giant component, we
start this exploration process  with $n-(1+O(\delta)){2d \over
d-2} \eps n$ tuples of which $n-(1+O(\delta)){2d(d-1)\over
d-2}\eps n$ are $d$-tuples and $(1+O(\delta))2d\eps n$ are
$(d-1)$-tuples. The number of vertices is therefore $dn -
(1+O(\delta)) {4d(d-1) \over d-2} \eps n$. In the notation of
Section \ref{explore} we have
$$ \Big | \NN_{0}^{(d)} - (n-{2d(d-1) \over d-2}\eps n) \Big | \leq \delta \eps n
\, , \qquad \qquad |\NN_{0}^{(d-1)} - 2d\eps n| \leq \delta \eps n
\, .$$ Fix
$$ T = (1+\eta)2(d-2)\eps^{-2}\log(\eps^3 n) \, .$$
As $|\NN_{t}^{(k)} - \NN_{t-1}^{(k)}| \leq 2$ for every $t$ and
$k$, and $T \leq \delta \eps n$ we learn from (\ref{uniform0})
that for all $t \leq T$ we have \begin{eqnarray} \label{add-d}
\prob_M \Big ( \eta_t \in \Net_{t-1}^{(d)} \mid \F_{t-1} \Big )
&\leq& {d (n - (1+O(\delta)){2d(d-1) \over d-2} \eps n ) \over dn
- (1+O(\delta)){4d(d-1) \over d-2} \eps n} \nonumber \\ &\leq& 1 -
(1+O(\delta))2(d-1)\eps \, ,\end{eqnarray}

\begin{eqnarray} \label{add-d-1} \prob_M \Big ( \eta_t \in \Net_{t-1}^{(d-1)} \mid \F_{t-1}
\Big ) &\leq& { (d-1)(1+O(\delta))2d\eps n \over dn -
(1+O(\delta)){4d(d-1) \over d-2}\eps n} \nonumber \\ &\leq&
(1+O(\delta)) 2(d-1)\eps \, .\end{eqnarray} By (\ref{incdef}) we
can bound $\prob_M(\xi_t = d-2 \mid \F_{t-1})$ above by
multiplying the right hand side of (\ref{add-d}) times $p$.
Similarly, we can bound $\prob_M(\xi_t = d-3 \mid \F_{t-1})$ above
by multiplying the right hand side of (\ref{add-d-1}) times $p$.
Therefore, we can stochastically bound from above $\xi_t$ by a
random variable $\beta$ taking the value $d-2$ with probability
${1 - (1+O(\delta))(2d-3)\eps \over d-1}$, the value $d-3$ with
probability $(1+O(\delta)) 2 \eps$ and otherwise the value $-1$.
Recall that $t_1$ denotes the first hitting time of $0$ by the
process $\{Y_t\}$. Lemma \ref{hitting2} then gives
$$ \prob_M \Big ( t_1 > T \Big ) \leq \eps
(n\eps^3)^{-(1+\eta)(1-O(\delta))(1-O(\eps))} \leq \eps
(n\eps^3)^{-(1+\eta/2)} \, ,$$ as long as $\delta$ is small
enough. Applying Lemma \ref{excessbd} and Lemma \ref{bottomline}
gives that
$$ \prob_M \Big ( |\C(v)| > (1+\eta){2(d-2) \over d-1} \eps^{-2}\log(\eps^3 n) \Big ) \leq  \eps
(n\eps^3)^{-(1+\eta/2)} \, ,$$ and as in the proof of Proposition
\ref{dregular} this yields that
$$ \prob_M \Big ( |\S_2| > (1+\eta){2(d-2) \over d-1} \eps^{-2}\log(\eps^3
n)\Big ) \leq (\eps^3 n)^{-\eta /2} \to 0 \, .$$

The proof that for every fixed $\ell \geq 2$
$$ \prob_M \Big ( |\S_\ell| < (1-\eta){2(d-2) \over d-1} \eps^{-2}\log(\eps^3
n)\Big ) \to 0 \, ,$$ goes by bounding the process $Y_t$ from
below by a process with independent increments. This is carried
out almost identically to the proof of Theorem \ref{lower} and we
omit the details. \qed \\

\section{The limiting distribution}\label{seclimit}

Recall the definitions of the processes $B^\lambda (\cdot)$ and
$W^\lambda(\cdot)$ in (\ref{process1}) and (\ref{process2}).
Throughout this section for a process $\{S_t\}$ indexed by
positive integers we write $S_t$ for $t \in \mathbb{R}$ to denote
the continuous linear interpolation of $S_t$.

%
Recall that $0=t_0 < t_1 < t_2 < \ldots$ are the times at which
$A_{t_j}=0$. Using the process $Y_t$ we define the process
$\widehat{Y}_t$ by $\widehat{Y}_0=Y_0 = d$ and for any $t \in
[t_j, t_{j+1})$ \be \label{fix}
 \widehat{Y}_t = \left \{
\begin{array}{ll}
Y_t, & \textrm{if } Y_t \geq Y_{t_j} \, , \\
Y_{t_j} & \textrm{otherwise} \, , \\
\end{array} \right .
\ee and $\widehat{Y}_t = \widehat{Y}_{dn/2}$ for any $t \geq
dn/2$. In this manner, the times $\{t_j\}$ are all the record
minima of the process $\{\widehat{Y}_t\}$. The main theorem of
this Section is the following:

\begin{thm}\label{convthm} Fix $\lambda \in {\mathbb R}$ and let $p={1+\lambda n^{-1/3} \over
d-1}$. Then as $n\to \infty$ we have that
$$n^{-1/3}\widehat{Y}_{((d-1)n^{2/3} \cdot)} \stackrel{d}{\Longrightarrow} (d-1)B^\lambda (\cdot)  \, , $$ where this convergence is on
finite intervals.
\end{thm}

Theorem \ref{convthm} states that $n^{-1/3}
\widehat{Y}_{((d-1)n^{2/3} \cdot)}$ converges to the process
$(d-1)B^\lambda$. It is thus natural to expect that ordered
excursions lengths of $\widehat{Y}_{((d-1)n^{2/3} \cdot)}$ above
past minima, will converge to excursions lengths of $B^\lambda$
above its past minima. Theorem \ref{criticaldist} essentially
follows from this assertion, but proving it requires some
technical work and we provide the details below. We postpone the
proof Theorem \ref{convthm} to the end of this section. \\

Fix some $s>0$ and let $C[0,s]$ be the space of continuous real
functions on $[0,s]$. Let $f \in C[0,s]$ and consider the set
$${\cal E} = \{(r,\ell)\subset [0,s] :
f(r)=f(\ell)=\min_{u \leq \ell}f(u) \textrm{ and } f(x)
> f(r) \quad \forall r < x < \ell\} \, .$$
This set defines excursions of $f$ above its past minima. To each
excursion $(r, \ell)$ we associate the length $\ell-r$. Since the
sum of excursion lengths is at most $s$, it is possible to order
them in a decreasing order $(\L_1, \L_2, \ldots)$. We call a point
$\ell$, such that $(r,\ell) \in {\cal E}$, an excursion {\em
ending} point. We say a function $f \in C[0,s]$ {\em good} if none
of its excursion ending points are local minima and if almost
every point in $[0,s]$ is contained in some excursion, i.e. for
almost every $x\in [0,s]$ there exists $(r, \ell) \in {\cal E}$
such that $r < x < \ell$. Given an integer $m$, consider the
function $\phi_m : C[0,s] \to {\mathbb R}^m$ defined by
$$ \phi_m (f) = (\L_1, \ldots , \L_m) \, .$$

\begin{prop}
If $f\in C[0,s]$ is good, then $\phi$ is continuous at $f$ with
respect to the $||\cdot||_{\infty}$ norm.
\end{prop}
\noindent {\bf Proof.} We prove for the case $m=1$. The proof for
$m>1$ is similar and we omit it. Let $f_n \in C[0,s]$ be a
sequence of functions such that $f_n \to f$. Consider the longest
excursion $(r, \ell)$ such that $\ell -r = \L_1 = \phi_1(f)$. As
for any $\eps>0$ small enough there exists $\delta>0$ such that
$f(x) > f(r) + \delta$ for $x \in (r + \eps, \ell -\eps)$ we
conclude that $\liminf _{n \to \infty} \phi_1(f_n) \geq \phi_1
(f)$. On the other hand, as almost every point in $[0,s]$ is
inside some excursion of $f$, for any $\eps >0$ we can find
excursions ending points $\ell_1, \ldots \ell_k$ of $f$ such that
$\ell_1 \leq \L_1 + \eps$, $s-\ell_k < \L_1 + \eps$ and $\ell_i -
\ell_{i-1} < \L_1 + \eps$ for $1 < i \leq k$. Since $f$ is good,
for any $\eps>0$ small enough we can find $\delta>0$ such that
there exists $x_i \in (\ell_i, \ell_i +\eps)$ such that $f(\ell_i)
- f(x_i) > \delta$ for all $i\leq k$. It follows that for large
enough $n$, the function $f_n$ has excursion ending points in the
intervals $(\ell_i, \ell_i + \eps)$. We conclude that $\limsup
_{n\to \infty} \phi_1(f_n) \leq \phi_1(f)$. \qed \\

%

\noindent {\bf Proof of Theorem \ref{criticaldist}.} See \cite{MP}
or \cite{RY} for general background on Brownian Motion and for the
proofs of the theorems we use in the following. Fix some $s>0$. It
is a classic fact that the zero set of Brownian motion has no
isolated points and is of $0$ measure with probability $1$. Also,
by a Theorem of Levy we know that $\{B(t) - \min _{y\leq t}B(y)
\}_t$ is distributed as $\{|B(t)|\}_t$, so we deduce that with
probability $1$, a Brownian motion sample path is good. By the
Cameron-Martin Theorem, with probability $1$ the process
$B^\lambda (\cdot)$ is good. As $\phi_m$ is continuous on almost
every sample point of $B^\lambda$, and $\phi_m ((d-1)B^\lambda) =
\phi_m (B^\lambda)$ we deduce by Theorem \ref{convthm} and Theorem
2.2.3 from \cite{D} that for any integer $m>0$
$$ ((d-1)n^{2/3})^{-1} \phi_m (\widehat{Y}_u) \stackrel{d}{\Longrightarrow} \phi_m (B^\lambda)\, .$$

In Section \ref{explore} we showed that the times $t_j$ are record
minima of $Y_t$. Hence, by (\ref{fix}), the lengths
$\{t_{j+1}-t_j\}$ are excursions lengths of $\widehat{Y}_u$ above
its past minima. Lemma \ref{excessbd} allows us to deduce
immediately that for any $s>0$,

\be\label{error} n^{-2/3} \Big | \Big \{ t \leq sn^{2/3} : \eta_t
\in \Act_{t-1} \hbox{ or } [\eta_t] \in \cup _{i=0}^{d-2}
\Net_{t-1}^{(i)} \Big \} \Big | \stackrel{d}{\Longrightarrow} 0 \,
.\ee Thus, if $t_{j+1} - t_j$ is the $\ell$-th largest excursion
ending before time $sn^{2/3}$, if $n^{-2/3}(t_{j+1}-t_j)$
converges in distribution to some random variable $\chi$, then
Lemma \ref{bottomline} and (\ref{error}) imply that the $\ell$-th
largest component completely explored before time $sn^{2/3}$,
normalized by $n^{-2/3}$, converges in distribution to ${\chi
\over d-1}$. As excursion lengths of $\widehat{Y}_{(n^{2/3}
\cdot)}$ are excursion lengths of $\widehat{Y}_{((d-1)n^{2/3}
\cdot)}$ times $(d-1)$ we learn that the sizes of components
discovered before time $sn^{2/3}$ in the exploration process,
normalized by $n^{-2/3}$ and ordered, converge in distribution to
the ordered excursion sizes of $W^\lambda[0,s]$.

We also need to handle the issue of the simplicity of the
resulting graph. The following lemma will be useful and is an
immediate consequence of Theorem $1$ and $2$ of \cite{BC}.

\begin{lemma} \label{enum} Let $d\geq 3$ and let $\bar{d}_1, \bar{d}_2 \in \{1,
\ldots,d \}^m$ be degree sequences of length $m$ such that each
sequence sum to an even number. Let $\prob_1$ be the distribution
of a uniform perfect matching on $\sum _{i=1}^m \bar{d}_1(i)$
vertices, divided to $m$ tuples such that the $i$-th tuple has
$\bar{d}_1(i)$ vertices in it. Similarly define $\prob_2$ using
degree sequence $\bar{d}_2$. Let $\Simp$ be the event that
contracting each tuple into a single vertex yields a simple graph.
Assume $d$ is fixed and $m \to \infty$. If $\bar{d}_1 = (d, \ldots
,d)$ and $\bar{d}_2$ has $(1-o(1))m$ entries with the value $d$
then
$$ \prob_2 (\Simp) = (1+o(1)) \prob_1 (\Simp) \, .$$
\end{lemma}


Fix a real number $s>0$ and consider a fixed interval $I \subset
[0,s]$. Let $\A_I$ denote the event $$ \A_I = \Big \{ n^{-2/3}
\Phi_m(\widehat{Y}_{(n^{2/3} \cdot)}) \in I \Big \} \, .$$
 For times $t<t'$ denote by $\S[t,t']$ the event that no loops
or parallel edges (either closed or open) were found between times
$t$ and $t'$ by the exploration process. The closed and open edges
inspected by the exploration process are a uniform random
matching, hence we have that $\prob (\S[0,dn/2]) = \prob(\Simp)$.
After $t$ steps of the exploration process the number of
$d$-tuples with $d$ neutral vertices is at least $n-2t$. Hence,
Lemma \ref{enum} shows that if $t=o(n)$ then $\prob ( \S[t,dn/2]
\mid \F_t ) = (1+o(1))\prob (\Simp)$. Thus, by conditioning on
$\F_{sn^{2/3}}$ we find that
$$ \prob \Big ( \A_I \,\cap \, \Simp \Big  ) = (1+o(1))\prob (\A_I) \prob (\Simp) \, .$$
Hence, when we condition on $\Simp$, component sizes discovered up
to time $sn^{2/3}$, normalized, also converge to excursions of
$W^\lambda[0,s]$. \\

Since we handled only components discovered before time $sn^{2/3}$
for some arbitrary large $s>0$, our final task for completing the
proof is to show that large components are typically found in the
beginning of the process, rather the end of it. The next lemma
completes the proof of the theorem. \qed \\

\begin{lemma} \label{tailcomp} Let $\C^{(sn^{2/3})}_1$ be the largest component which we
started exploring {\em after} time $sn^{2/3}$. Then for any
$\alpha > 0$ we have

\be \label{finito} \lim _{s \to \infty} \limsup_{n \to \infty}
\prob \Big ( | \C^{(sn^{2/3})}_1 | \geq \alpha n^{2/3} \Big ) = 0
\, .\ee
\end{lemma}

\noindent{\bf Proof of Lemma \ref{tailcomp}.} Let $\hat{t}_0 >
sn^{2/3}$ be the first time larger than $sn^{2/3}$ at which
$A_{\hat{t}_0} = 0$ and let $m = \sum _{k=1}^d
\NN_{\hat{t}_0}^{(k)}$. We continue the exploration process on a
graph that has $m$ tuples, of varying sizes between $1$ and $d$,
in which the number of $k$-tuples is $\NN_{\hat{t}_0}^{(k)}$.
After finishing the exploration process we again contract each
tuple to a vertex to form the graph $G^*_m$ on the vertex set $U$
of cardinality $m$ . The components discovered before $\hat{t}_0$
together with $G^*_m$ form $G^*(n,d,p)$. Our analysis will show
that from any starting vertex $u \in U$, the drift of the process
$\{Y_t\}$ is too small to have components of size $\alpha
n^{2/3}$.

Fix some small $\delta>0$ and denote by $\A$ the event
$$ \A = \Big \{  \NN_t^{(d)} \leq n-(1-3\delta)t \, , \quad \NN_t^{(0)} \geq (1-3\delta)pt  \quad \forall t  \in [sn^{2/3}, \delta n] \Big \} \, .$$
By (\ref{incdef}) and (\ref{uniform0}) we have that for all $t$
\begin{eqnarray*}
\E [ \xi_t \mid \F_{t-1} ] &\leq& p { \sum _{k=2}^d k(k-1)
\NN_{t-1}^{(d)} \over dn -2t + 1} -1 \\
&\leq& p \Big [ {d(d-1)\NN_{t-1}^{(d)} + (d-1)(d-2)(n-
\NN_{t-1}^{(d)} - \NN_{t-1}^{(0)}) \over dn-2t+1} \Big ] - 1 \, ,
\end{eqnarray*}
where the last inequality is due to the fact that $\sum
_{k=2}^{d-1} \NN_{t-1}^{(k)} \leq (n-\NN_{t-1}^{(d)} -
\NN_{t-1}^{(0)})$. We now substitute $p={1+\lambda n^{-1/3} \over
d-1}$ and condition on $\A$. A straightforward calculation gives
that for $t  \in (sn^{2/3}, \delta n]$, we have
$$ \E [ \xi_t \mid \F_{t-1}, \A ] \leq \Big (1+\lambda n^{-1/3}\Big )\Big (1-{[(d-2) + O(\delta)] t \over d(d-1) n}\Big ) -1 \, .$$
We deduce that if $s=s(\delta, \lambda)>0$ is large enough, then
for all $t \in (sn^{2/3}, \delta n]$,

\be \label{condbound} \E [ \xi_t \mid \F_{t-1}, \A ] \leq -
\delta^{-1} n^{-1/3} \, .\ee

Assume we start exploring at time $\hat{t}_0+1$ the tuple of a
vertex $u\in U$ (i.e., $w_{\hat{t}_0+1}$ is in the tuple
corresponding to $u$). Denote by $\C(u)$ the connected component
of $u$ and by $\gamma$ the stopping time
$$ \gamma = \min \{ t > 0 : Y_{\hat{t}_0 + t} = Y_{\hat{t}_0} - N(w_{\hat{t}_0 +1})\} \,
.$$ By bounding $U_j \leq t_j - t_{j-1}$ and $V_j\leq t_j -
t_{j-1}$ in Lemma \ref{bottomline} we get

$$ |\C(u)| \leq \Big ( {2\over d-1} + 1 \Big ) \gamma + 1 \, .$$
By optional stopping and (\ref{condbound}), since $N(w_{\hat{t}_0
+1})\leq d$, we have that $\E[\gamma \wedge \delta n \mid \A] \leq
\delta d n^{1/3}$ as long as $s$ is large enough. By
(\ref{deviation2}) and (\ref{deviation3}) of Lemma \ref{extradev}
we have that for $n$ large enough $\prob ( \A ^c ) \leq n^{-1}$.
Also, part $1$ of Theorem \ref{criticalupper} implies that
$\prob(\gamma > \delta n) \leq n^{-1}$ for large enough $n$.
Hence,
\begin{eqnarray*} \E \gamma &\leq& dn \prob (\gamma > \delta n) + \E [
\gamma{\bf 1}_{\{\gamma < \delta n\}} ] \leq d + \E [\gamma \wedge
\delta n ] \\ &\leq& d + \E [\gamma \wedge \delta n \mid \A ] +
\delta n \prob (\A^c) \leq (d+1)\delta n^{1/3} \, ,
\end{eqnarray*}
for large enough $s>0$. The same analysis works for any $u \in U$
and so we learn that that $\E |\C(u)| \leq O(\delta) n^{1/3}$ for
all $u \in U$ . Thus for any fixed $\alpha>0$ we have
$$ \prob ( |\C(u)| > \alpha n^{2/3} ) \leq O(\delta) n^{-1/3} \, ,$$
where the constants in the O-notation depend on $\alpha$ and $d$.

Let $X$ be the random variable counting the number of $u \in U$
such that $|\C(u)|>\alpha n^{2/3}$. As $m \leq n$ we have proved
that $\E X \leq O(\delta)n^{2/3}$. Observe that
$|\C^{(sn^{2/3})}_1|
> \alpha n^{2/3}$ implies that $X>\alpha n^{2/3}$. Hence
$$ \prob \Big ( |\C^{(sn^{2/3})}_1| > \alpha n^{2/3} \Big ) \leq
O(\delta) \, .$$ Since $\delta>0$ was
arbitrary and $s$ was large enough depending only on $\delta$ and $\lambda$, this concludes our proof. \qed \\

We now turn to the proof of Theorem \ref{convthm}. For the proof
we use a standard functional central limit theorem for martingales
(see \cite{D}, Theorem 7.2):

\begin{thm} \label{FCLT}
Let
$$ \{ X _{m,k}, {\cal F}_{m,k} : 1 \leq k \leq m \} \, ,$$
be a martingale difference array. For any $\ell \leq m$ let
$$ V_{m,\ell} = \sum _{k=1}^\ell \E [ X_{m,k}^2 \mid {\cal F}_{m,k-1} ] \,
$$
be the quadratic variation process, and
$$ Z_{m,\ell} = \sum _{k=1}^\ell X_{m,k} \, .$$ If
\begin{enumerate}
\item $|X_{m,k}| \leq \delta_m$ with $\delta_m \to 0$, and \item
for each $t\in [0,1]$ we have $V_{m,\lfloor mt \rfloor} \to t$ in
probability as $m\to \infty$,
\end{enumerate}
then $Z_{m,(mt)} \stackrel{d}{\Longrightarrow} B(t),$ where
$B(\cdot)$ is standard Brownian motion, and $Z_{m,(\cdot)}$ is the
continuous linear interpolation of $Z_{m,k}$.
\end{thm}

\noindent {\bf Proof of Theorem \ref{convthm}.} Since
$|n^{-1/3}Y_u - n^{-1/3}\widehat{Y}_u| \leq dn^{-1/3}$ for all
$u\geq 0$, it suffices to prove the convergence for the process
$\{Y_u\}$. Fix some $s>0$, take $m=m_n=\lfloor sn^{2/3} \rfloor$
and denote $\xi_k^* = \E [\xi_k \mid \F_{k-1}]$. Consider the
martingale difference array,
$$ X_{m,k} = m^{-1/2}( \xi_k - \xi_k^*), \quad k \leq m \, .$$
We have that for any $\ell \leq m$,
\begin{equation}\label{intheway} Z_{m,\ell} = \sum _{k=1}^\ell
X_{m,k} = m^{-1/2} Y_\ell - m^{-1/2}\sum _{k=1}^\ell \xi_k^* \,
.\end{equation}

As $|X_{m,k}|=O(n^{-1/3})$, condition $1$ of Theorem \ref{FCLT} is
satisfied. Putting $\eps = \lambda n^{-1/3}$ in $(i)$ and $(ii)$
of  Corollary \ref{drift} gives that
$$ \sup _{k \leq m_n} | \xi_k^* | \to 0 \quad \hbox{in $L_1$ and in probability.}$$
Hence by $(iii)$ of Corollary \ref{drift} we get
$$ \sup _{k \leq m_n} \E [ (\xi_k - \xi_k^*)^2 \mid \F_{k-1} ] \to (d-2) \quad \hbox{ in probability,}$$
as $n \to \infty$. Thus for any $t \in [0,1]$ we have
$$ m^{-1} \sum_{k=1}^{\lfloor mt \rfloor} \E[
(\xi_k - \xi_k^*])^2 \mid \F_{k-1} ] \to (d-2)t  \quad \hbox{ in
probability.} \,  .$$ In the notation of Theorem \ref{FCLT}, it
follows that $V_{m,\lfloor mt \rfloor} \to (d-2)t$ in probability.
We conclude by Theorem \ref{FCLT} that
$$ Z_{m,(mt)} \stackrel{d}{\Longrightarrow} B((d-2)t) \,
. $$ An immediate computation with Part (ii) of Corollary
\ref{drift} shows that
$$ \E \sum _{k=1}^m |\xi_k^* - \E \xi_k| = O(1) \, ,$$ which by
the triangle inequality gives that

\be \label{1storderap} m^{-1/2} \E \max _{k_0 \leq m} \Big |
\sum_{k=1}^{k_0} \xi_k^* - \sum _{k=1}^{k_0} \E\xi_k \Big |=
O(n^{-1/3}) \, .\ee Part (i) of Corollary \ref{drift} with $\eps
=\lambda n^{-1/3}$ implies that for any $t \in [0,1]$,
$$ m^{-1/2} \sum _{i=0}^{tm} \E \xi_i \longrightarrow \lambda
t \sqrt{s} - {(d-2) t^2 s^{3/2} \over 2d(d-1)} \, .$$ We conclude
by (\ref{1storderap}) that

\be \label{usethis} \prob \Big ( \sup _{t\in [0,1]} \Big |
m^{-1/2} \sum_{k=1}^{\lfloor tm \rfloor} \xi_k^* - \lambda
\sqrt{s} t + {(d-2) s^{3/2}t^2 \over 2d(d-1)} \Big | > n^{-1/6}
\Big ) \longrightarrow 0 \, .\ee Rearranging (\ref{intheway})
using (\ref{usethis}) gives that for any fixed $s>0$
$$ {n^{-1/3} \over \sqrt{s}}Y _{(n^{2/3} ts)}
\stackrel{d}{\Longrightarrow} B((d-2)t) + \lambda t \sqrt{s} -
{(d-2) t^2 s^{3/2} \over 2d(d-1)} \, .$$ Multiplying by $\sqrt{s}$
and using Brownian scaling gives
$$ {n^{-1/3}}Y _{(n^{2/3} ts)}
\stackrel{d}{\Longrightarrow} B( (d-2) ts) + \lambda ts - {(d-2)
(ts)^2 \over 2d(d-1)} \, .$$ By Brownian scaling and the
definition of $B^\lambda$ we deduce that
$$ {n^{-1/3}}Y _{((d-1) n^{2/3}\cdot )}
\stackrel{d}{\Longrightarrow} (d-1) B^\lambda (\cdot) \, ,$$
which concludes our proof. \qed \\

\section{Concluding Remarks}\label{secconc}

\begin{itemize}

\item It is natural to ask whether the bounds in Proposition
\ref{dregular} are tight. In light of Theorem \ref{criticalupper}
we would expect that for $\lambda \in \mathbb{R}$ there exists a
constant $c=c(\lambda)$ such that for {\em any} $d$-regular graph
$G$ and $A>0$ we have
$$ \prob \Big ( | \C_1(G_p) | > An^{2/3} \Big ) \leq e^{-cA^3} \, .$$
The authors currently know how to prove this for some particular
cases, for instance, expander graphs.

\item It is an interesting topic for further research to find a
quenched version of Theorem \ref{criticaldist}. Recall that
$|\gamma_1|$ is the longest excursion above past minima of the
process $B^\lambda$ defined in (\ref{process1}). Let $D(n,d)$
denote the number of simple $d$-regular graphs on $n$ vertices and
set $p={1 + \lambda n^{-1/3} \over d-1}$. We expect that for small
$\eps_1>0, \eps_2>0$ and any $s>0$ and $n$ large enough at least
$(1-\eps_1)D(n,d)$ of the $d$-regular graphs $G$ on $n$ vertices
satisfy
$$ \Big | \prob ( |\C_1 (G_p)| < sn^{2/3} ) - \prob ( |\gamma_1|
\leq s ) \Big | \leq \eps_2 \, .$$

\item Assume now $d=d(n)$ grows with $n$. We proved that when
$d(n)$ is a fixed constant, then $G(n,d(n),p)$ is mean field
around ${1 \over d(n)-1}$. The same result holds for $d(n)=n-1$
since this is just the usual $G(n,p)$ model. It seems plausible
that for all such sequences (assuming $nd(n)$ is even) the same
conclusion still holds.

\end{itemize}

\section*{Acknowledgments}
We are grateful to Itai Benjamini for posing this problem and for
many insightful conversations. Part of the work was done while the
first author was an intern at Microsoft Research.

We also thank Jian Ding, Nathan Levy, Alex Smith and Nick Wormald
for useful comments and suggestions.

\bigskip \noindent
{\bf Asaf Nachmias}: \texttt{asafnach(at)math.berkeley.edu} \\
Department of Mathematics\\
UC Berkeley\\
Berkeley, CA 94720, USA

\medskip \noindent
{\bf Yuval Peres}: \texttt{peres(at)stat.berkeley.edu} \\
Microsoft Research,
One Microsoft way,\\
Redmond, WA 98052-6399, USA.

\end{document}